\documentclass[11pt]{amsart}
\usepackage{geometry}                % See geometry.pdf to learn the layout options. There are lots.
\usepackage{graphicx}
\usepackage{amssymb}
\usepackage{epstopdf}
\usepackage{comment}

\usepackage{framed}
\newcommand{\bfra}{\begin{framed}\noindent}
\newcommand{\efra}{\end{framed}}

\usepackage[T1]{fontenc}
\usepackage{pxfonts}

\usepackage{color}

\newcommand{\red }{\color{red}}

\newcommand{\lab}{\label}

\newcommand{\bconj}{\begin{conjecture}}
\newcommand{\econj}{\end{conjecture}}
\newcommand{\bques}{\begin{question}}
\newcommand{\eques}{\end{question}}

\newcommand{\ben}{\begin{enumerate}}
\newcommand{\een}{\end{enumerate}}

\newcommand{\bea}{\begin{eqnarray}}
\newcommand{\ba}{\begin{array}}
\newcommand{\bean}{\begin{eqnarray*}}
\newcommand{\ea}{\end{array}}
\newcommand{\eea}{\end{eqnarray}}
\newcommand{\eean}{\end{eqnarray*}}
\newcommand{\beq}{\begin{equation}}
\newcommand{\eeq}{\end{equation}}
\newcommand{\bthm}{\begin{thm}}
\newcommand{\ethm}{\end{thm}}
\newcommand{\blem}{\begin{lem}}
\newcommand{\elem}{\end{lem}}
\newcommand{\bprop}{\begin{prop}}
\newcommand{\eprop}{\end{prop}}
\newcommand{\bcor}{\begin{cor}}
\newcommand{\ecor}{\end{cor}}
\newcommand{\bdfn}{\begin{dfn}}
\newcommand{\edfn}{\end{dfn}}
\newcommand{\brem}{\begin{rem}}
\newcommand{\erem}{\end{rem}}
\newcommand{\bpf}{\begin{proof}}
\newcommand{\epf}{\end{proof}}
\newcommand{\bfact}{\begin{fact}}
\newcommand{\efact}{\end{fact}}
\newcommand{\bobs}{\begin{obs}}
\newcommand{\eobs}{\end{obs}}
\newcommand{\bexam}{\begin{exam}}
\newcommand{\eexam}{\end{exam}}
\newcommand{\bclaim}{\begin{claim}}
\newcommand{\eclaim}{\end{claim}}

\newtheorem{thm}{Theorem}
\newtheorem{prop}[thm]{Proposition}
\newtheorem{lem}[thm]{Lemma}

\newtheorem{cor}[thm]{Corollary}
\newtheorem{dfn}[thm]{Definition}
\newtheorem{rem}[thm]{Remark}
\newtheorem{fact}[thm]{Fact}
\newtheorem{claim}[thm]{Claim}
\newtheorem{obs}[thm]{Observation}
\newtheorem{exam}[thm]{Example}

\newtheorem{conjecture}[thm]{Conjecture}
\newtheorem{question}[thm]{Question}

\newtheorem*{condition'}{Condition 2'}

 \newtheoremstyle{claimstyle}%
   {}%             space above
   {}%             space below
   {\normalfont}%     body font
   {}%                indent
   {\itshape}%        header font
   {.}%               punctuation
   { }%               space after head
   {\thmnote{#3}}%    typeset note only.

\theoremstyle{claimstyle}

\alph{enumii} \roman{enumiii}

             \def\cB{\mathcal B}       
\def\cH{\mathcal H}             \def\cF{\mathcal F}       
           \def\cM{\mathcal M}        
                       \newcommand{\J}{\mathcal{J}}
              \def\cS{\mathcal S}             
                    
\def\cD{\mathcal D}

                \def\Z{{\mathbb Z}}      \def\R{{\mathbb R}}
\def\C{{\mathbb C}}                      \def\oc{{\hat \C}}
    			\def\D{{\mathbb D}}  
 					\def\H{{\mathbb H}}  

\newcommand{\amsc}{{\mathbb C}}

\newcommand{\cbar}{\hat{{\mathbb C}} }

\def\a{\alpha}                \def\b{\beta}             \def\d{\delta}
                          
                           \def\l{\lambda} 
\def\La{\Lambda}                         \def\Om{\Omega}
               \def\sg{\sigma}
\def\Th{\Theta}                          
\def\ka{\kappa}

\newcommand{\lam}{\lambda}
\newcommand{\ep}{\varepsilon}
\newcommand{\ph}{\varphi}
\newcommand{\al}{\alpha}
\newcommand{\ga}{\gamma}

\def\lt{\left}                \def\rt{\right}

\def\ben{\begin{enumerate}}
\def\een{\end{enumerate}}

\def\${$\displaystyle}
\def\ph{\varphi}
\def\sp{\medskip}   
\def\fr{\noindent}
 \def\HD{\text{{\rm HD}}} 
          \def\P{\text{{\rm P}}}

\def\HypDim{{\rm HD_{hyp}}}

\def\escape{\mathcal I(f)}

\def\expd{\gamma}
\def\tstar{\Theta}

\newcommand{\jul}{J(f)}
\newcommand{\fat}{F(f)}
\newcommand{\sing}{S(f)}
\newcommand{\conical}{J_r^{sph}(f)}

\newcommand{\ehyp}{{\textsc E}--hyperbolic }
\newcommand{\ehypty}{{\textsc E}--hyperbolicity }
\newcommand{\eh}{{\textsc E}--}

\newcommand{\exc}{{\mathcal E}_f}

\newcommand{\post}{{\mathcal P}(f)}

\newcommand{\pft}{{\mathcal{L}}_t}

\newcommand{\Tract}{\Om}

\def\1{1\!\!1}

\def\den{\rho}

\newcommand{\pf}{{\mathcal{L}}}
\newcommand{\npf}{\mathcal{\hat L}}

\def\rad{J_r(f)}

\title[Thermodynamic Formalism for Transcendental Functions]{Thermodynamic Formalism and Geometric Applications \\ for \\ Transcendental Meromorphic and Entire Functions}

\author{Volker Mayer}
\address{Volker Mayer, Universit\'e de Lille, D\'epartement de Math\'ematiques, UMR 8524 du CNRS, 59655 Villeneuve d'Ascq Cedex, France}
\email{volker.mayer@univ-lille.fr \newline \hspace*{0.42cm}}

\thanks{
Research of the second named author supported in part by the
Simons Grant: 581668.}

\author{Mariusz Urba\'nski}
\address{Mariusz Urba\'nski, Department of Mathematics, University of North Texas, Denton, TX 76203-1430, USA}
\email{urbanski@unt.edu \newline \hspace*{0.42cm} \it Web: \rm www.math.unt.edu/$\sim$urbanski}

\begin{document}
\maketitle
%\section{}
%\subsection{}

\tableofcontents
%---------------------------------------------------------------------------------------------------------------------------------------------------%

\pagebreak 

\section{Introduction}

\

Originating from statistical physics, the dynamical theory of thermodynamic formalism was brought to mathematics,
particularly to study expanding and hyperbolic dynamical systems, 
primarily by  R. Bowen \cite{Bow75}, D. Ruelle \cite{Rue78}, Ya. Sinai \cite{ Si72}, and P. Walters \cite {W1} in the 1970's. This theory provides an excellent framework for probabilistic description of the chaotic part of the dynamics and, in the context of smooth (particularly conformal) expanding/hyperbolic dynamical systems, gives a rich and detailed information about the geometry of expanding repellers, limit sets of Kleinian groups and iterated function systems, and Julia  sets of holomorphic dynamical systems.  More precisely, by establishing the existence and uniqueness of Gibbs and equilibrium states, and
studying spectral and asymptotic properties of corresponding
Perron-Frobenius operators, it permits to show that dynamical systems are "strongly" mixing (K-mixing, weak Bernoulli), have exponential decay of correlations, satisfy the Invariant Principle Almost Surely, in particular satisfy the Central Limit Theorem, and the Law of Iterated Logarithm.
Furthermore, by studying the topological pressure function of geometric potentials, particularly its regularity properties (real analyticity, convexity), this theory
gives a precise information about the fractal geometry of Julia and limit sets. Particularly, R. Bowen initially showed in \cite{Bow79} that the Hausdorff dimension of the limit set of a co--compact quasi--Fuchsian group is given by the unique zero of the appropriate pressure function. His result and its numerous versions commonly bear the name of Bowen's Formula ever since. Bowen's results easily carry through to the case of expanding (hyperbolic) rational functions providing a closed formula for the Hausdorff dimension of their Julia sets.
D. Ruelle, positively answering a conjecture of D. Sullivan (see \cite{Su4}--\cite{Su3}), proved in \cite{Rue82} that this dimension depends in a real analytic way on the function. 

For hyperbolic, and even much further beyond, rational functions, and more general distance expanding maps, the theory of thermodynamic formalism is now well developed and established, and its systematic account can be found in \cite{PUbook}
(see also \cite{Wal82, Zins2000}, \cite{MRU}, \cite{KU20}).
The present text concerns transcendental entire and meromorphic functions. For these classes of functions
many differences and new phenomena pop up that do not occur in the case of rational maps. 
The following two properties of transcendental functions show from the outset that the outlook of these classes is indeed totally different than the one of rational functions.
\ben
%\item[$\bullet$] 
\item[-]  Whereas the singularities of hyperbolic rational maps stay away form their Julia sets,
for transcendental functions one always has to deal with the singularity at infinity.
\item[-] Transcendental functions have infinite degree.
\een
One immediate consequence of the later fact is that for transcendental functions there is no measure of maximal entropy, which is one of the
central objects in the theory of 
rational functions. Particularly for polynomials, where this measure coincides with harmonic measure viewing from infinity, and also 
for endomorphisms of higher dimensional projective spaces.
Another consequence is that all Perron--Frobenius, or transfer, operators of a transcendental meromorphic functions are always defined
by an infinite series. This is the reason that, even for such classical functions as exponential ones $f_\l(z)=\l e^z$, this operator
taken in its most natural sense, is not even well--defined.

K. Bara\'nski first managed to overcome these difficulties and presented a
thermodynamical formalism for the tangent family in \cite{Baranski95}.
Expanding the ideas from \cite{Baranski95} led to \cite{ku1}, where Walters expanding
maps and Bara\'nski maps were introduced and studied. One important feature of the maps
treated in \cite{Baranski95} and  \cite{ku1} was that all analytic inverse branches were
well-defined at all points of Julia sets. This property dramatically fails
for example for entire functions as $f_\l(z)=\l e^z$ (there are no
well-defined inverse branches at infinity). To remedy this situation, the periodicity of $f_\l$ was exploited to project
the dynamics of these functions
down to the cylinder and the appropriate thermodynamical formalism was developed
 in \cite{UZ03} and \cite{UZ04}. This approach has been adopted to other periodic
 transcendental functions; besides the papers cited above, see also
 \cite{CoiSko07a, CoiSko07b, KU04, MyUrb05.2, UZ07, DUZ07} and the survey \cite{KU08}.
 
\sp
The first general theory of thermodynamic formalism for transcendental meromorphic and entire functions was laid down in the year 2008 in \cite{MyUrb08}. 
 \cite{MUmemoirs} containes a complete treatment of this approach.
 It handled all the periodic functions cited above in a uniform way and
  went much farther beyond.   
 The most important key point in these two papers was to replace the standard Euclidean metric by an appropriate Riemannian metric. Then the power series defining the Perron--Frobenius operators of geometric potentials becomes comparable to the Borel series and can be controlled by means of Nevanlinna's value distribution theory.
 
 \sp For a large class of transcendental entire functions whose set of singularities is bounded, quite an optimal approach to thermodynamic formalism was laid down and developed in \cite{MUpreprint3}. 
 
 It was observed in this paper that, for these entire functions, the transfer operator entirely depends on the geometry of the logarithmic tracts,
 in fact on the behavior of the boundary of the tracts near infinity.  The best way to deal with the often fractal behavior of the tracts near infinity was 
 by adapting the concept of integral means, a classical and powerful tool in the theory of conformal mappings.

\sp This text provides an overview of the (geometric) thermodynamic formalism for transcendental meromorphic and entire functions 
  with particular emphasis  on geometric/fractal aspects such as Bowen's Formula expressing the hyperbolic dimension as a unique zero of a pressure function and the behavior of the latter when the transcendental functions vary in an analytic family. 
 
There are some several important and interesting topics closely related to the subject matter of our exposition that will nevertheless not be treated at all or will be merely briefly mentioned in our survey. For example, this exposition only briefly indicates that thermodynamic formalism has been successfully developed for random transcendental dynamical systems; see \cite{MyUrb2014, Mayer2016aa}, comp. also \cite{UZ18} for non--hyperbolic random dynamics of transcendental functions. 
Non--hyperbolic functions will not be in the focus of our current exposition either but we would like to bring reader's attention to some relevant papers that include \cite{UZ07}, \cite{MyUrb10.2} and \cite{UZ18}. Discussing all these topics at length and detail would increase the length of our survey substantially, making it too long, and would lead us too far beyond of what we intended to focus on in the current survey.

\sp
We would like to thank the referee for his valuable remarks which influenced the final version of the paper.

\section{Notation}

Frequently we have to replace Euclidean metric by some other 
 Riemannian metric  $ d\sigma =\gamma\,|dz|$. A natural choice is the
 spherical metric in which case the
 density with respect to Euclidean metric is $\ga (z) =1/(1+|z|^2)$. More generally, we consider metrics of the form
\beq \label{2.8 full}
d\sigma (z) =d\sigma_\tau (z) = \frac{|dz|}{1+|z|^\tau} \quad , \quad \tau \geq 0.
\eeq
They vary between euclidean and spherical metrics when $\tau\in [0,2]$.
If such a metric is used only away from the origin, then one can use the simpler form 
\beq \label{2.8}
d\tau (z) = |z|^{-\tau}|dz|.
\eeq
We denote by $D_\sg (z,r)$ the open disk with center $z$ and radius $r$ with respect metric $\sg$.
If $\sg$ is the spherical metric then this disk is also denoted by $D_{sph} (z,r)$ and for
the standard euclidean metric $\D (z,r)$. We also denote 
$$
\D_R =\D(0,R)
$$
and $$\D^* _R = \C\setminus \overline \D (z,R).$$
The symbol
$$
A(r,R):=\D_R\setminus \overline \D_r
$$
is used to denote the annulus centered at $0$ with the inner radius $r$ and the outer radius $R$. 

The derivative of a function $f$ with respect to a  
 Riemannian metric  $ d\sigma =\gamma
\,|dz|$ is given by 
\beq\label{derivative R metric}
|f'(z)|_\sigma =\frac{d\sigma (f(z))}{d\sigma (z)} =|f'(z)|
\frac{\gamma (f(z))}{\gamma (z)}.
\eeq
When the metric $\sg$ has the form \eqref{2.8 full} or \eqref{2.8}
then $d\sigma$ only depends on $\tau$ and we will identify $\sg$ and $\tau$ and
write $|f'(z)|_\tau$ instead of $|f'(z)|_\sigma $. 
Therefore,

$$|f'(z)|_\tau = \frac{|f'(z)|}{1+|f(z)|^\tau}(1+ |z|^\tau)
\quad\text{and}\quad
|f'(z)|_\tau = \frac{|f'(z)|}{|f(z)|^\tau} |z|^\tau
$$
in the case of the simpler form \eqref{2.8}.
When $\tau =2$ then we also write $|f'(z)|_{sph}$.

\sp Besides this, we use common notation such as $\C$ and $\cbar$ for the Euclidean plane and the Riemann sphere respectively. Another common notation is 
$$
A \asymp B.
$$
As usually, it means that the ratio $A/B$ is bounded below and above by strictly positive and finite constants
that do not depend on the parameters involved. The corresponding inequalities up to a multiplicative constant are 
denoted by 
$$
A\preceq B
\  \  \ {\rm and} \  \  
A\succeq B.
$$
Also,
$$
\text{dist} (E,F)
$$ denotes the Euclidean distance between the sets $E,F \subset \C$.

%---------------------------------------------------------------------------------------------------------------------------------------------------%
\section{Transcendental functions, hyperbolicity and expansion}

%\subsection{Generalities on transcendental functions}
We consider transcendental entire or meromorphic functions. Such a function  $f:\C\to \cbar$ can have two types of singularities:
asymptotical and critical values. We refer to \cite{BE08.1} for the classification of the different types of singularities, known as Iversen's classification,
denote by $S(f)$ the closure of the set of critical values and finite asymptotic values of $f$.

Transcendental functions are very general and one is led, actually forced, to consider reasonable subclasses.
The class $\cB$ of bounded type functions consists of all meromorphic functions for which the set $S(f)$ is bounded. 
Bounded type entire functions have been introduced and studied in \cite{EL92},
$\cB$ is also called the  Eremenko--Lyubich class. It contains an important subclass, called Speiser class, which consists of all meromorphic functions for which the set $S(f)$ is finite. 

\subsection{Dynamical preliminaries}

\label{section2} For a general introduction of the dynamical aspects of meromorphic functions
we refer to the survey article of Bergweiler \cite{Bergweiler-survey} and the book \cite{KU20}. We
collect here some of its properties, primarily the ones we will need in the sequel. The Fatou
set of a meromorphic function $f:\amsc\to \cbar$ is denoted by
$\fat$. It is defined as usually to be the set of all points $z\in \C$
for which there exists a neighborhood $U$ of $z$ on which all the iterates
$f^k$, $k\geq 1$, of the function $f$ are defined and form a normal family.
The complement of this set is the
Julia set ${\hat {\mathcal J}}(f)=\oc \setminus \fat$. We write
$$
\jul = {\hat {\mathcal J}}(f)\cap \C.
$$
\index{$\jul$}
By Picard's theorem, there are
at most two points  $\xi\in \cbar$ that have finite backward orbit
${\mathcal O}^- (\xi)=\bigcup_{n\geq 0}f^{-n}(z_0)$. The set of
these points is the exceptional set $\exc$. In contrast to the
case of rational maps it may happen that $\exc \subset {\hat {\mathcal J}}(f)$. Iversen's theorem \cite{Iv1914, Nevbook70} asserts that every point $\xi\in \exc$
is an asymptotic value. Consequently, $\exc$ is contained in $\sing$. the set
of critical and finite asymptotic values of $f$.
\index{$\sing$} \index{$\post$}
The post-critical set
$\post$ is defined to be the closure in the complex plane $\C$ of
$$
\bigcup_{n\geq 0} f^n \big(\sing \setminus f^{-n}(\infty ) \big)\; .
$$
%{\green
This set can contain the whole Julia set.
\bdfn\label{tame}
If $\jul \setminus \post \neq \emptyset$ then $f$ is called tame.
\edfn
%}%end gree,

The Julia set contains several dynamically important subsets. First, there is
the \emph{escaping set}
$$
\escape = \{z\in \C \, ; \; \text{$f^n(z)$ is defined for all $n$ and } \lim_{n\to\infty} f^n(z) =\infty \}.
$$
This set is not always a subset of the Julia set, it may contain Baker domains.
However, for entire functions of bounded type $\escape \subset \jul$ (\cite[Theorem 1]{EL92}).
More important for us is the following set.

\bdfn\label{defi conical euclidean} The
\it radial (or conical) \rm Julia set $\rad$ of $f$ is the set of
 points $z\in \jul$ such that there exist
$\delta >0$ and an unbounded sequence $(n_j)_{j=1}^\infty$ of positive integers such that the sequence $\big(|f^{n_j}(z)|\big)_{j=1}^\infty$ is bounded above and the map
$$ f^{n_j} : U_j\longrightarrow    \D (f^{n_j} (z) , \delta )$$
is conformal,
where $U_j$ is the connected component of $f^{-n_j} (\D (f^{n_j} (z) , \delta ))$ containing $z$.
\edfn

There are other definitions of radial sets in the literature. While the present definition is in the spirit of the one from \cite{SU14},
the radial points  in \cite{Rempe09-hyp} are defined by means of spherical disks. Namely,  $z\in \conical$ if $z\in \jul$
 if  there exist
$\delta >0$ and an unbounded sequence $(n_j)_{j=1}^\infty$ of positive integers such that
$$ f^{n_j} : U_j\longrightarrow    D _{sph} (f^{n_j} (z) , \delta )$$
is conformal
where $U_j$ is the connected component of $f^{-n_j} (D_{sph} (f^{n_j} (z) , \delta ))$ containing $z$.
Right from these definitions it is easy to see that $\rad\subset\conical$. Also,
\beq\label{inclusion rad escape}
\rad \subset \jul\setminus \escape\, .
\eeq
The differences between all of these radial sets are dynamically insignificant in the sense that they all have the same Hausdorff dimension and this dimension coincides with the hyperbolic dimension, which we define right now.

\bdfn \label{defi hyp dim}
The hyperbolic dimension of a meromorphic function $f:\C\to\cbar$, denoted by $\HypDim (f)$,  is 
$$
\HypDim (f) = \sup _K \HD (K)
$$
where
the supremum is taken over all hyperbolic sets $K\subset\C$, i.e. over all compact sets $K\subset \C$ such that $f(K)\subset K$ and $f|_K$ is expanding.
\edfn

\blem
$\HypDim (f)= \HD (\rad ) =\HD (\conical)$.
\elem

\bpf
Let $K$  be a hyperbolic set. Then, following \cite[Section 5]{PUbook} especially Lemma 5.1.1, there exists $\eta >0$
 such that 
 $$
 f_{|\D (z, \eta )} \text{ injective and } \quad f(\D(z, \eta ) ) \supset \D (f(z) , \eta ) 
 \quad \text{for all } \, z\in K.
 $$
This shows that $K\subset \rad$ and thus $\HypDim (f) \leq \HD (\rad )$.
Since $\rad\subset\conical$ we also have $\HD (\rad ) \leq \HD (\conical )$. The conclusion comes now from
 the result in  \cite{Rempe09-hyp} which says that $\HD (\conical ) = \HypDim (f)$.
\epf

%\smallskip

\subsection{Hyperbolicity and expansion}
There are several notions of hyperbolic transcendental functions in the literature
(see for example \cite{Zhen15}). The following definition is used fairly frequently.

\bdfn\label{defi hyp}
A meromorphic function $f:\C\to\cbar$ is called hyperbolic if and only if
\beq\label{hyper bdd type}
\post \text{ is bounded and} \;  \; \post \cap \jul =\emptyset.
\eeq
\edfn

Notice that then $f\in \cB$, i.e. it is of bounded type. The following notion
has been used by G. Stallard in \cite{St1990}. Later it was considered in  \cite{MyUrb08, MUmemoirs} where the, somehow misleading, name topologically hyperbolic was used. Since it is based on the euclidean distance, let us call it \emph{Euclidean hyperbolic} here.

\bdfn\label{defi hyp topo} A meromorphic function $f$ is called
\it \ehyp \index{\ehyp} \rm if
$$
\text{dist} (\jul ,\post ) >0\,.
$$
\edfn

\fr 
Clearly, a hyperbolic function is \ehyp  but the later notion is much more general.
For example, the function $f(z)=2-\log 2 +2z -e^z$
is \ehyp and has a Baker domain (see \cite{BwInvariant95}). Other examples arise naturally in the context of Newton maps. This has been observed in \cite{BFJK2019} and  $f(z)=z-\tan z$ is a different exemple of \ehyp  function that is not hyperbolic.

%radial
\smallskip

For \ehyp  functions every non-escaping point of the Julia set is a radial point.
Together with \eqref{inclusion rad escape} it follows that in this case
we have equality between these type of points:
\beq \label{eq radial entire}
\rad = \jul \setminus \escape .
\eeq

%end radial

For rational functions, \ehypty is equivalent to the property of being expanding.

\bdfn\label{defi expa}
A meromorphic function $f:\C\to\cbar$ is called \it expanding  \rm if and only if there are two constants $c>0$ and $\expd >1$ such that
$$
|(f^n)'(z)|\ge c\expd ^n \quad for \; all  \;\; z\in \jul \setminus
f^{-n}(\infty  )\; .
$$
%A  \ehyp and expanding function is called \it hyperbolic.
\edfn

\smallskip

The function $f(z)=z-\tan z$ is  \ehyp  and $|f'(z)|\to 1$ as $\Im f(z) \to \infty$.
Since there are vertical lines in $\jul$, it follows that this function is not expanding. Thus,
contrary to the case of rational maps, \ehypty and
expanding are not equivalent for transcendental functions. It is shown in  \cite{St1990} that
every entire \ehyp 
function $f$ satisfies $\lim_{n\to \infty}|(f^n)'(z)|\to \infty$ for all $z\in \J_f$ and under some conditions
the expanding property follows from
\ehypty (see Proposition 4.4 in \cite{MUmemoirs}).

\bexam\label{exam 8}
Let $0<c<1/e^3$. Then the Fatou function $f(z)=z-\log c +e^{-z}$ is not hyperbolic 
but it is \ehyp  and expanding.

\rm
In order to verify this statement, we recall the classical argument that $f$ is semi-conjugate via $w=e^{-z}$ to the map
$g(w)=cwe^{-w}$ (see for example \cite{Rippon08}). By the choice of  the constant  $0<c<1/e^3$, the origin is an attracting fixed point of the map $g$ and
a simple estimation allows to check that 
\beq\label{eq d_3}
g(\overline\D_3)\subset \D_3.\eeq Consequently, the half space
$\{ \Re z \geq -\log 3\}$ is contained in a Baker domain of $f$ and the Julia set $\jul\subset \{ \Re z < -\log 3\}$. 
Now, a simple  estimate shows that
$$
|f'(z)|\geq 2 \  \  \  \text{for all} \  \ 
  z  \  \  {\rm with} \  \  \Re z \leq -\log 3.
$$
Consequently $f$ is expanding on its Julia set.

It remains to check that $f$ is \ehyp . The function $f$ has no finite asymptotic value and its
 critical points are
 $c_k=2\pi i k$, $k\in \Z$. It follows from
\eqref{eq d_3} that there exists $\rho> -\log 3$ such that $\Re f^n(c_k) \geq \rho$ for every $n\geq 0$ and  $k\in \Z$.
This shows that $f$ is indeed \ehyp.
\eexam

% Disjoint type
\subsection{Disjoint type entire functions}

%\sp\fr 
For entire functions there is a relevant strong form of hyperbolicity called \emph{disjoint type}, the notion that
first implicitly appeared in \cite{B07} and then was explicitly studied in several papers including \cite{RRRS, Rempe09}.
Disjoint type functions are of bounded type. So, let $f\in \cB$ be an entire function and let $R>0$ such that 
$\sing \subset \D _R$. Up to normalization we can assume that $R=1$. Then $f^{-1}(\D^*)$ consists of countably many mutually disjoint unbounded Jordan domains $\Omega_j$ with real analytic boundaries such that $f:\Om_j\to \D^*$ is a covering map (see \cite{ EL92}). In terms of the classification of singularities, this means that $f$ has only \emph{logarithmic singularities} over infinity.
These connected components of $f^{-1}(\D^*)$ are called \emph {tracts} and 
the restriction of $f$ to any of these tracts $\Om_j$ has the special form 
\beq\label{11}
\text{$f_{|{\Om_j}} =\exp\circ{\tau_j}\;\;$
where $\;\;\ph_j=\tau_j^{-1} :\cH=\big\{z\in\C:\Re (z) >0\big\} \longrightarrow \Om_j$}
\eeq
is a conformal map. Later on we often assume that $f$ has only finitely many tracts:
\beq\label{8}
f^{-1}(\D^*)=\bigcup_{j=1}^N \Om_j\,.
\eeq
Notice that this is always the case if the function $f$ has finite order. Indeed, if $f$ has finite order then the Denjoy-Carleman-Ahlfors Theorem (see \cite[p. 313]{Nevbook74}) states that $f$ can have only finitely many direct singularities and so, in particular, only finitely many logarithmic singularities over infinity. 

\bdfn\label{defi disjoint type}
If $f\in \cB$ is entire such that 
\beq\label{1}
S(f)\subset \D \quad  \  \text{and} \quad  \  \bigcup_j \overline\Om_j \cap \overline \D  =\emptyset   \;\; , \;  \text{  equivalently} \  f^{-1}(\overline{\D^*})=\bigcup_j \overline\Om_j \subset  \D^*, 
\eeq
then $f$ is called a \emph{disjoint type} function. 
\edfn

This definition is not the original one but it is consistent with the disjoint type models in Bishop's paper \cite{Bishop-EL-2015}.
The function $f$ is then indeed of disjoint type in the sense of \cite{B07, RRRS, Rempe09}. 
It is well known that for every $f\in \cB$  the function $\l f$, $\l\in \C^*$, is of disjoint type provided $\l$ is small enough (see 
\cite{BK-2007} and \cite[p.261]{Rempe09}).
Also, the Julia set of a disjoint type entire function is a subset of its tracts and therefore only the restriction of the function to these
tracts is relevant for the study of dynamics of such a function near the Julia set.

\smallskip

Besides functions of class $\cS$ and $\cB$ we consider the following subclass of bounded type entire functions called class $\cD$.
In this definition,  
\beq\label{rectangles} Q_T=\big\{0<\Re z<4T \; ; \; -4T<\Im z < 4T\big\}\quad ,\quad T >0.\eeq

 \bdfn\label{9}
An entire function $f:\C\to \C$ belongs to class  $\cD$ if it is of disjoint type, has only finitely many tracts (see  \eqref{8})
and if, for every tract, the function $\ph$ of \eqref{11} satisfies
\beq\label{(4.2)}
|\ph(\xi )|\leq M |\ph (\xi ')| \quad \text{for all \ $\;\xi , \xi ' \in Q_T\setminus Q_{T/8}$}, 
\eeq
 for some constant $M\in (0,+\infty)$ and every $T\geq 1$.  \edfn

%\sp
%-----------------------------------------------------------------------------------------------

\section{Topological Pressure and Conformal Measures}\label{sec Pressure and conf}

This section is devoted to two crucial objects: the topological pressure and conformal measures.
Compared to the case of rational functions, they both behave totally differently in the context of 
transcendental functions. For example, since transcendental functions have infinite degree, the topological pressure evaluated at zero is always infinite. Also, the existence
of the pressure and, even more importantly, of conformal measures is not known in full generality for meromorphic functions.

\subsection{Topological Pressure} A standard argument, based on mixing properties (see Lemma 5.8 in \cite{MUmemoirs}), shows that for a \ehyp  meromorphic function $f:\C\to \C$
the following number, which might be finite or infinite,
\beq \label{marre 2}
\P_{\tau} (t) := \limsup_{n\to \infty} \frac 1n\log \sum_{f^n(z)=w} |(f^n)'(z)|_{\tau}^{-t}
\eeq
does not depend on the point $w\in\jul$. However, this number may depend on the metric $\tau$ and it clearly depends on the parameter $t>0$.

\bdfn\label{topo press}
Let $f:\C\to \C$ be \ehyp  meromorphic function. The topological pressure of $f$ evaluated at $t>0$ with respect to
the metric  $\tau$ as defined in  \eqref{2.8 full} is the (possibly infinite) number $\P_{\tau} (t)$ defined by formula \eqref{marre 2}. When $\tau=2$, i.e. $d\sg_\tau$ is the spherical metric, then we also write $\P_{\tau} (t)=\P_{sph} (t)$.
\edfn

Given a meromorphic function $f:\C\to \C$, a number $\tau\ge 0$, and a parameter $t\ge 0$, we say that the topological pressure $\P_{\tau} (t)$ exists if the number defined by formula \eqref{marre 2} is independent of $w$ for some "sufficiently large'' set of points $w\in \jul$.

The most general result on the existence of topological pressure going beyond \ehyp  functions is due to Bara\'nski, Karpi\'nska and Zdunik \cite{BKZ-Bowen}.
They work with spherical metric and call a meromorphic function $f$ \emph{exceptional} if and only if it has a (Picard) exceptional value $a$
in the Julia set and $f$ has a non--logarithmic singularity over $a$.

\bthm[\cite{BKZ-Bowen}]\label{thm sph pressure}
Let $f:\C\to \C$ be either a meromorphic function in class $\cS$ or a non--exceptional and tame function in class $\cB$. Then the limit
$$
\P_{sph} (t) := \lim_{n\to \infty} \frac 1n\log \sum_{f^n(z)=w} |(f^n)'(z)|_{sph}^{-t}
$$
exists (possibly equal to infinity) for all $t>0$ and does not depend on $w$ where $w$ is a good pressure starting point $w\in \C$ whose precise meaning  is given in \cite[Section 4]{BKZ-Bowen}. 

If $f$ is tame then every $w\in \jul\setminus \post$ is such a good point. Also, if $f\in \cB$ is \ehyp, then each point $w\in \jul$ is good.
\ethm

It is also shown in \cite{BKZ-Bowen} that the pressure function has the usual natural properties. 

\bprop [\cite{BKZ-Bowen}]
Under the assumptions of Theorem \ref{thm sph pressure},
$\P_{sph} (0) =+\infty$ and $\P_{sph} (2)\leq 0$, and thus 
\beq\label{Th pressure}
\Th_{sph}:=\inf \{t>0 \; : \; \P_{sph} (t) < \infty\}\in [0,2]
\eeq
In addition:
$$
\P_{sph} (t)=+\infty
\  \  \ {\rm for \ all } \  \  \  t<\Th_{sph} \  \ {\rm and } \  \  \ \P_{sph} (t)<+\infty \  \  \ {\rm for \ all } \  \  \  t>\Th_{sph}. 
$$
The resulting function 
$$
(\Th_{sph},\infty)\ni  t\longmapsto \P_{sph}(t)
$$ 
is non--increasing and convex, hence continuous.
\eprop

Notice that this result does not provide any information
about the behavior of the pressure function at the critical value $t=\Th_{sph}$. For classical families, such as the exponential family, the pressure  at 
$\Th_{sph}$ is infinite. Curious examples of functions that behave differently at the critical value are provided in \cite{MyZd}. We will come back to such examples later in Theorem \ref{thm myzd1}.

%\medskip

\subsection{Conformal measures and transfer operator} 
Conformal measures were first defined and introduced by Samuel Patterson in  his seminal paper \cite{Pat1} (see also \cite{Pat2}) in the context of Fuchsian groups. Dennis Sullivan extended this concept to all Kleinian groups in \cite{su82, Su2, Su3}. He then, in the papers \cite{Su4, Su5, Su6}, defined conformal measures for all rational functions of the Riemann sphere $\oc$. He also proved their existence therein. Both Patterson
and Sullivan came up with conformal measures in order to get an understanding of geometric measures, i.e. Hausdorff and packing ones. Although already Sullivan noticed that there are conformal measures for Kleinian groups that are not equal, nor even equivalent, to any Hausdorff and packing (generalized) measure, the main purpose to deal with them is still to understand Hausdorff and packing measures but goes beyond. 

Conformal measures, in the sense of Sullivan have been studied in greater detail in \cite{DU2}, where, in particular, the structure of the set of their exponents was examined. We do this for our class of transcendental functions.

Since then conformal measures in the context of rational functions have been studied in numerous research works. We list here only very few of them appearing in the early stages of the development of their theory: \cite{DU LMS}, \cite{DU3}, \cite{DU4}. Subsequently the concept of conformal measures, in the sense of Sullivan, has been extended to countable alphabet iterated functions systems in \cite{MU1} and to conformal graph directed Markov systems in \cite{MauUrb03}. It was furthermore extended to transcendental meromorphic dynamics in \cite{ku1}, \cite{UZ03}, and \cite{MyUrb08}.
See also \cite{MUmemoirs}, \cite{UZ04}, and \cite{BKZ18}.
Lastly, the concept of conformal measures found its place also in random dynamics; we cite only \cite{MSUspringer}, \cite{MyUrb2014}, and \cite{UZ18}.

\bdfn\label{4.3.3}
Let $f:\C\to \hat \C$ be a meromorphic function. A Borel probability measure $m_t$ on $\jul$ is called $\l |f'|_\tau^{t }$--conformal if 
$$ 
m_t (f(E)) = \int _E \l |f'|_\tau^{t } dm_t .
$$
for every Borel $E\subset \jul $ such that the restriction $f_{|E}$ is injective.
%If $\l=1$ then $m_t$ is called geometric conformal, or geometric $t$--conformal, measure.
The scalar $\l$ is called the conformal factor and, if $\l=1$, then $m_t$ is called a $t$--conformal measure.
\edfn
  
If $f$ has  a $\l |f'|_\tau^{t }$--conformal measure $m_t$ and if $f$ is \ehyp  then, using Koebe's Distortion Theorem,
we get for all $w\in \jul$ that
\beq\label{marre 1}
1\geq 
\sum_{z\in f^{-1}(w)} m_t(U_z) 
\asymp \l ^{-1}m_t( \D(w,r)) \sum _{z\in f^{-1}(w)}  |f'(z)|_\tau ^{-t}  ,
\eeq
where for every $z\in f^{-1}(w)$, $U_z$ is the connected component of $f^{-1}(\D(w,r))$ containing $z$. Consequently, the 
series on the right hand side of \eqref{marre 1} is well defined. This allows us to introduce the corresponding transfer, or Perron--Frobenius--Ruelle, operator.
Its standard definition as an operator acting on the space  $ C_b (\jul)$  of continuous bounded functions on the Julia set $\jul$ is the following.

\bdfn\label{defi rpf operator}
Let $f:\C\to \hat\C$ be a \ehyp  meromorphic function. Fix $\tau>0$ and $t>0$. The transfer operator of $f$ with (geometric) potential $\psi:=-t\log |f'(z)|_\tau$, $t>0$, is defined by
\beq\label{13/3 1}
\pft g (w):=\sum_{f(z)=w}e^{\psi (z)}g(z) =\sum_{f(z)=w}|f'(z)|_\tau ^{-t} g(z) \quad , \quad w\in \jul \; , \; g\in C_b (\jul).
\eeq
\edfn

Note that \eqref{marre 1} does not imply boundedness of the linear operator $\pft$. This crucial issue will be discussed in the next section.
Let us simply mention here that, iterating the inequality \eqref{marre 1}
(which is possible if we assume $f$ to be \ehyp) shows that we have the following relation between
 the pressure and the conformal factor  $\l$:
 \beq\label{2/4 1}
\P_\tau (t)\leq \log \l \, .
\eeq

On the other hand, if the transfer operator, in fact its adjoint operator  $\pft^*$, is well defined,
then $m_t$ being a conformal measure equivalently means that
$m_t$ is an eigenmeasure of $\pft^*$ with eigenvalue $\l$: 
$$
\pft^* m_t = \l m_t.
$$
%If in addition the Ruelle--Perron--Frobenius Theorem \ref{theo pfr} holds then 
%the eigenvalue $\l=\l_t$ is determined by the pressure $\P(t)$ since then we have equality in \eqref{2/4 1}. We will see later that this is not always the case.

%The name \emph{geometric} conformal measure is justified by the fact that the case $\l=1$ is
%most interesting for the study of the geometry of Julia sets.

As defined, the measure $m_t=m_{\tau, t}$ 
does depend on the metric $\tau$. Given $\tau ' \neq \tau$ and corresponding Riemannian metrics (see \eqref{2.8}), we then have
\beq\label{23/3 9} 
\frac{ dm_{\tau ', t}}{ dm_{\tau, t}}(z) = \frac{|z|^{(\tau -\tau')t}}{ \int _{\jul}  |\xi|^{(\tau -\tau')t} d m_{\tau, t}(\xi)}  
\eeq
provided the above integral is finite. For example, this allows one to get spherical conformal measures
as soon as we have conformal measures $m_{\tau, t}$  for a $\tau$--metric with $\tau \leq 2$; this will be the case later in the results Theorem \ref{theo pfr} and Theorem \ref{21/3 1}. Indeed, the formula
\beq\label{23/3 90} 
\frac{ dm_{sph, t}}{ dm_{\tau, t}}(z) := \frac{|z|^{(\tau -\tau')t}}{ \int _{\jul}  |\xi|^{(\tau -\tau')t} d m_{\tau, t}(\xi)}  
\eeq
defines a spherical conformal probability measure $m_{sph, t}$. 
But then it may happen that the corresponding density (Radon--Nikodym derivative) $d\mu_t / dm_{sph ,t}$ in 
Theorem \ref{theo pfr} or Theorem \ref{21/3 1}
is no longer a bounded function.

\smallskip
\subsection{Existence of conformal measures}

As we have already said, for rational functions, Denis Sullivan proved  in \cite{su82} that every rational function
admits a conformal measure with conformal factor $\l=1$.
For transcendental functions
this is not so in full generality and this is again because of 
the singularity at infinity. 
In general, a conformal measure is obtained
by a (weak) limit procedure and one has to make sure that the mass does not
escape to infinity when passing to the limit. 

\sp 

There are two particular cases where natural $t$--conformal measures do exist.
First of all, there are different types of meromorphic functions for which 
 the  (normalized) spherical Lebesgue measure is a $2$--conformal measure.
 This is the case for functions $f$ with $\jul=\cbar$ and for those having a Julia set
 of positive area such as the functions of the sine family. This is a result of Curtis McMullen \cite{McM87}; we will come back to it and to its generalizations in greater detail in Section \ref{2 conf meas}. 
 
 The other particular case is formed by meromorphic functions having as their Julia sets the real line $\R$ or a geometric circle, and thus having a natural $1$--conformal measure. Functions of this type arise among inner functions studied by  Aaronson \cite{Aa97}
and Doering--Man\'e \cite{DM89}.

\sp Coming now to the general case, the relation between topological pressure and the existence of conformal measures
has been studied in \cite{BKZ18}. The hypotheses of this paper 
are again those of Theorem \ref{thm sph pressure} and thus it
goes far beyond (\eh) hyperbolic functions.
 Theorem C of that paper contains the following general statement for the existence of $t$--conformal measures.

\bthm[\cite{BKZ18}]
Let $f:\C\to\oc$ be either a meromorphic function in class $\cS$ or a non--exceptional and tame function in class $\cB$.
If $\P_{sph} (t)=0$ for some $t>0$, then $f$ has a $|f'|^t_2$--conformal measure,
i.e. a $t$--conformal measure, with respect to the spherical metric.
\ethm

For \ehyp and expanding function there exists a general construction of conformal measures. It allows us to produce conformal measures, defined with respect to adapted $\tau$--metrics,  with various conformal factors $\l$. The proof of Proposition 8.7 in \cite{MUpreprint3} along with Section 5.3 in \cite{MUmemoirs} yield the following.

\bthm\label{t120200730}
Let $f:\C\to\oc$ be  \ehyp and expanding. Assume that $t>0$ and $\tau$ are  such that 
$$
\| \pft \1 \|_\infty < +\infty \quad \text{ and }  \quad \lim_{|w|\to \infty \, ,\; w\in \jul} \pft \1 (w) =0 .
$$
Then there exist a $\l |f'|_\tau ^t$--conformal measure with $\l = e^{\P_\tau (t)}$.
\ethm

The first hypothesis of this theorem tells us that we have a ``good'' well defined bounded linear transfer operator. The second hypothesis can be used to prove tightness of an appropriate sequence of purely atomic measure, which in turn allows us to produce, as its weak* limit, a desired conformal measure. Then Theorem~\ref{t120200730} follows.

\smallskip
\subsection{Conformal measures on the radial set and recurrence}

For rational functions the behavior of conformal measures on the radial set is fairly well understood.
For example, it has been studied in \cite{DMNU} and in \cite[Section 5]{McM-HdimII}, and most of the arguments from these papers can be adapted to the transcendental case.

\bthm \label{thm 27...} Let $m_t$ be a $\l|f'|_\tau^t$--conformal measure of a meromorphic function $f:\C\to\oc$ such that $m_t(\rad )>0$.
Then
$$
\text{$m_t(\rad ) =1\;\; $ ,}
$$
 $m_t$ is ergodic, $m_t$ almost every point has a dense orbit in $\jul$ and $m_t$ is a unique $\l|f'|_\tau^t$--conformal measure. More precisely, if $m$ is a $\rho |f'|_\tau^t$--conformal measure %with $m(\rad )>0$, 
 then $\l =\rho$ and $m=m_t$.
%thus, conformal mesures are dissipative or recurrent? ($\tau \leq 2$ for distortion and then one can always pass to spherical?)
\ethm

\bpf
The radial Julia set has been defined in Definition \ref{defi conical euclidean}. 
For any $z\in \rad$, let $\d (z)>0$ be the number $\d$
and let $(n_j)_{j\geq1}$ be the sequence associated to $z$, both according to Definition \ref{defi conical euclidean}. Define then
$$
J_r (f , \d ):=\big\{ z\in \rad \; : \;\; \d (z)\geq 2\d \text{ and } \sup_{j\ge 1} \big\{|f^{n_j}(z)|\big\}\leq 1/\d \big\}.
$$
Then
\beq\lab{120200407}
\rad =  \bigcup _{\d >0} J_r (f , \d )
\eeq
and, if $m_t (\rad ) >0$, then $m_t (J_r(f,\d ) ) >0$ for some $\d>0$.

For all $z\in \rad$, consider the blow up mappings 
$$
f^{n_j}:V_j(z)\longrightarrow \D \big(f^{n_j}(z),2\d\big) \; , \;\; j\ge 1,
$$
where $V_j(z)$ is the connected component of $f^{-n_j}\big(\D \big(f^{n_j}(z) , 2\d\big)\big)$ containing $z$. Let
$$
U_j(z):=V_j(z)\cap f^{-n_j}\big(\D \big(f^{n_j}(z) ,\d\big)\big).
$$
Then Koebe's Distortion Theorem applies for the map $f^{n_j}$ on $
U_j(z)$. In fact, what we need is a bounded distortion for the
derivatives taken with respect to the Riemannian metric $d\tau$. This however is a straightforward consequence of Koebe's Theorem (see \cite[Section 4.2]{MUmemoirs}).
Therefore
$$
m_t(U_j(z))\asymp \l ^{-n_j}|(f^{n_j})'(z)|^{-t}_{\tau} m_t (\D (f^{n_j}(z), \d )).
$$
Now, since conformal measures are positive on all non--empty open sets relative to $\jul$, we conclude that for every $\d >0$ there exists a constant $c>0$ such that
$$ 
m_t( \D(w,\d ))\geq c 
$$
for every $w\in \jul\cap B(0,1/\d )$.
This shows that 
\beq\label{2/4 2} 
m_t(U_j(z))\asymp\l ^{-n_j}|(f^{n_j})'(z)|^{-t}_{\tau}
\eeq
for every $z\in J_r(f,\d )$ and every $j\geq 1$
with comparability constants depending on $\d $ only.

Having this estimate we now can proceed exactly as in  \cite[Theorem 5.1]{McM-HdimII}.
If $\nu_t$ is any  $\eta |f'|_\tau^t$--conformal measure then \eqref{2/4 2} also holds with $\nu_t, \eta$
instead of $m_t, \rho$, and with other appropriated constants depending on $\d $ only. Hence, for every $z\in J_r(f,\d)$,
$$
\frac{ m_t( U_j(z))}{ \nu_t( U_j(z))} \asymp 1 \quad \text{for every $j\geq 1$}.
$$
Since in addition $\lim_{j\to\infty}{\rm diam}(U_j (z))=0$ and all $U_j(z)$, $z\in J_r(f,\d)$, $j\ge 1$, have shapes of not ``too much'' distorted balls, we conclude that the measures  $m_t$ and $\nu_t$ are equivalent (mutually absolutely continuous) on $J_r (f , \d )$. Invoking, \eqref{120200407}, we deduce that these two measures are equivalent on $\rad$. This is not the end of the proof yet but the interested reader is referred to the original proof in \cite{McM-HdimII}.
\epf

\smallskip

Recall that the Poincar\'e's Recurrence Theorem asserts that, given  $T:X\to X$ measurable dynamical system preserving a finite measure, for every measurable set $F\subset X$ and almost every point $x\in F$, the point $T^n(x)$ is in $F$ for infinitely many $n\ge 1$. A conformal measure $m$ is called \emph{recurrent} if the conclusion of the Poincar\'e recurrence theorem holds for it. 
In the case where the Perron--Frobenius--Ruelle theorem holds then,
due to the existence of probability invariant measures, commonly called Gibbs states, equivalent to the conformal measure, the later is always recurrent. 
By Halmos' Theorem \cite{Hal47}, recurrence is equivalent to \emph{conservativity} which means that
there does not exist a measurable wandering set of positive measure, i.e. a measurable set $W$ with $m(W)>0$
and such that 
$$
f^{-n} (W)\cap f^{-m} (W)=\emptyset
$$ 
for all $n>m\geq 0$.

\bthm\label{thm recurrence} Assume that the transcendental function $f:C\to\oc$ has $m_t$, a $\l|f'|_\tau^t$--conformal measure. Then
\ben
\item[-] $m_t$ is recurrent and this holds if and only if $m_t(\rad)=1$ or
\item[-] $m_t$--almost every point is in $\escape$ or its orbit is attracted by $\post$.
\een
\ethm

%\smallskip

\bpf
If $m_t ( \jul \setminus \post)=0$ then the second conclusion holds. So we may assume from now on
that $m_t ( \jul \setminus \post)>0$. Notice that then $f$ is tame and thus there exist  $D=\D(w,r)$, a disk centered at some point $w\in\jul$ and such that 
$$
\D(w,2r)\cap\post=\emptyset.
$$
Assume that there exists $W\subset \jul \setminus \post$  a wandering set of positive measure.  
Since all omitted values are in $\post$, there exists $N$ such that 
$$
W'=f^N (D)\cap W
$$ 
is a wandering set of positive measure. But then $W''= f^{-N}(W') \cap D$ is a wandering set of positive measure 
contained in $D$. Conformality, bounded distortion, and the fact the $W''$ is wandering, give 
$$
1\geq \sum_{n\geq 0} m_t(f^{-n}(W'')) \asymp m_t (W'') \sum_{n\geq 0} \pft^n \1 (z)
\asymp \frac{m_t (W'')}{m_t (D)} \sum_{n\geq 0} m_t(f^{-n}(D)) \, .
$$
The series in the middle is what is usually called the Poincar\'e series and we see that it is convergent for the exponent $t$. 
Now, a standard application of the Borel--Cantelli Lemma shows that a.e. $z$ is in at most
finitely many sets $f^{-n}(D)$ or, equivalently, only for finitely many $n$ we have $f^n(z)\in D$.
Since this true for every such disk $D$, it follows that 
\beq\label{eq cvb 1}
z\in \escape \quad \text{or} \quad f^n(z) \to \post \quad \text{for $m_t$ a.e. $z\in \jul$.}
\eeq
This also shows that $m_t(\rad )=0$ in this case since, as we have seen in Theorem \ref{thm 27...}, if 
$m_t(\rad )>0$ then $m_t$ a.e. orbit has a dense orbit in $J_f$ which contradicts \eqref{eq cvb 1}
since $f$ is a tame function.

\smallskip

The other possibility is that $\jul \setminus \post$ does not contain  a wandering set of positive measure.
Then $m_t$ is conservative hence recurrent on $\jul \setminus \post$. 
Let 
$$
V_\ep (A) := \{z\in \C \; : \; dist (z , A) \leq\ep\},
$$ 
$$
V^c_\ep (A):=\jul \setminus V_\ep (A)
$$ 
and consider the open set
$$U_\ep = \D(0, 1/\ep ) \cap V^c_\ep (\post ) \quad , \quad \ep >0.$$
If $\ep>0$ is small enough, $U_\ep\cap \jul \neq \emptyset$, and then $m_t(U_\ep)>0$.
On the other hand, recurrence implies that 
$$
m_t (U_\ep) = m_t \big(\{z\in U_\ep \, : \; f^n(z)\in U_\ep \text{ for infinitely many $n$'s}\}\big).
$$
The set of points $z$ such that, for some $\ep >0$, $z\in U_\ep$ and  $f^n(z)\in U_\ep$
for infinitely many $n$'s is a subset of $\rad$.  Therefore, $m_t(\rad) \geq m_t(U_\ep) >0$
and then, by Theorem \ref{thm 27...}, $m_t(\rad )=1$.
Notice also that then $m_t$ is recurrent on the whole Julia set since $m_t(\jul\setminus \rad ) =0$. 
\epf

%\smallskip

Every rational function has a $t$--conformal measure of minimal exponent $t=\d_p$; see \cite{DU2}. Shishikura \cite{Shi98}
gave the first examples of some polynomials $p$ for which this exponent is maximal, i.e. $\d_p=2$. For them the corresponding conformal measures are not recurrent. Up to our best knowledge it is unknown whether
there exist polynomials, even rational functions, $p$ with $\d_p<2$ and with non--recurrent $\d_p$--conformal measures. However, there are such 
quadratic like examples; see Avila--Lyubich \cite{AL08}, and the first globally defined, i.e. on the whole complex plane, (transcendental meromorphic) functions having
such behavior were produced in \cite[Theorem 1.4]{MyZd}. Notice that these examples are even hyperbolic
and their number $\Th$ (see Theorem \ref{thm myzd1}) is equal to the minimal exponent $\d_p$.

\bthm[\cite{MyZd}]\label{theo main' MyZd} 
There exist disjoint type entire functions $f:\C\to\oc$  
of finite order, with $\Th\in (1,2)$, that do not have any recurrent $\Th$--conformal measure with conformal factor $\l=1$. 
\ethm

\fr 
In fact  \cite[Theorem 1.4]{MyZd} states that these functions do not have $\Th$--conformal measures supported
on the radial Julia set. But this is equivalent to non--recurrence by Theorem \ref{thm recurrence}.

\smallskip

\subsection{$2$-conformal measures} \label{2 conf meas}
We finally discuss the special case of $2$--conformal measures.
As already mentioned above, for many transcendental, especially entire, functions
 the spherical Lebesgue measure $m_{sph}$ of the Julia set is positive and thus it is a natural $2$--conformal measure.
 In this case each of the following possibilities can occur:
 \begin{itemize}
 \item[-] $m_{sph}$ is recurrent and $m_{sph}(\rad ) =1$.
 \item[-] $m_{sph}(\escape \cap \jul )=1$.
 \item[-] $m_{sph}(\escape \cap \jul )=0$ and $f^n(z)\to \post$ for $m_{sph}$--a.e. $z\in \C$.
 \end{itemize}
 
 \smallskip
 
 Let us first discuss the recurrent case for which the postcritically finite map $f(z)=2\pi i e^z$
 is a typical example having the property that $m_{sph}(\rad ) =1$.
 For this function, the Julia set is the whole plane. From
  a classical zooming and Lebesgue density argument (see for example the proof of \cite[Theorem 8]{EL92})
  follows that this always holds provided that the radial set is positively charged.
  
\bprop
Let $f:\C\to\oc$ be a meromorphic function.
If $m_{sph}$ is a $2$--conformal measure and if $m_{sph} (\rad  )>0$, then $\jul = \C\,.$
\eprop

\fr 
For the third possibility, i.e. where the escaping set is not charged but where a.e. orbit is attracted by the post--critical set,
we have some results due to  Eremenko--Lyubich \cite[Section 7]{EL92}.

\bthm[\cite{EL92}]\label{criteria EL}
Let $f\in\cB$ be an entire function of finite order having a finite logarithmic singular value.
Then $m_{sph} (\escape ) =0 $ and there exists $M>0$ such that 
$$
\liminf_{n\to\infty} |f^n (z)| <M \quad \text{for a.e.  $z\in \C$.}
$$
\ethm

Given this result combined with Theorem \ref{thm recurrence} we see that there are several possibilities. Assume that $f$ satisfies the hypotheses of 
Theorem \ref{criteria EL} and that the Julia set of $f$ has positive area. Then, either the spherical Lebesgue measure is supported on
the radial Julia set or a.e. orbit is attracted by the post--critical set. 

As typical examples we can consider again the exponential family. As already mentioned,  $f(z)=2\pi i e^z$ is a recurrent example.
Totally different is $f(z)=e^z$.  Misiurewicz showed in \cite{Mi81} that $\jul=\C$ and Lyubich proved in \cite{Lyu87} that this function is not ergodic.
Consequently $m_{sph} (\rad )=0$ and thus a.e. orbit is attracted by the orbit of $0$, the only finite singular value.

\sp Plenty of entire functions have the property $m_{sph}(\escape \cap \jul )=1$ (and $\fat\neq \emptyset$).
Initially, McMullen showed in \cite{McM87} that the Julia set of every function from the sine--family $\a \sin (z)+\b$, $\a\ne 0$,
has positive area. This result has been generalized in many ways and to many types of entire functions;
 see \cite{AB12, BC16, Si15, Berg18}. The authors of these papers did not really deal with Julia but with the escaping set and, as a matter of fact, they showed that
\beq\label{area} 
area (\escape \cap \jul ) >0 \,.
\eeq
Since the escaping set is invariant, it suffices now to normalize properly the spherical Lebesgue measure restricted to $\escape$
in order to get the required $2$--conformal measure that is entirely supported on the escaping set.

\sp

%\vfill
%\pagebreak
%---------------------------------------------------------------------------------------------------------------------------------------------------%
\section{Perron--Frobenus--Ruelle Theorem, Spectral gap and applications}\label{sec transf operator}

The whole thermodynamic formalism relies on the transfer operator and its properties.
We recall that this operator has been introduced in Definition \ref{defi rpf operator}.
In fact, this definition treats only the most relevant \emph{geometric potentials}.
More general potentials $\psi$ were considered in  \cite{MUmemoirs}. They
are obtained as a sum of a geometric potential plus an additional H\"older function. This class of potentials 
has its importance for the multifractal analysis (see the Chapters 8 and 9 of  \cite{MUmemoirs}) of conformal measures and their invariant versions.
 In the present text we restrict ourselves to geometric potentials, so to functions of the form
 \beq\label{potential coboundary}
\psi := -t\log |f'| + b -b\circ f
\eeq
for some appropriate function $b:\jul \to \R$ (or $\C$). This coboundary is crucial since it allows us
to deal with different Riemannian metrics on $\oc$. We start by investigating elementary examples to make this transparent.

\smallskip

We already have mentioned in the introduction that the ``naive'' transfer operator is not always well defined. Let us consider the simplest entire function $f(z)=\l e^z$
and a potential $\psi := -t\log |f'|$ without coboundary. Then, for all $w\neq 0$ and parameter $t$,
%\beq\label{expo infty}
$$
\pft\1(w)
=\sum_{z\in f^{-1}(w)}|f'(z)|^{-t}
=\sum_{z\in f^{-1}(w)}|w|^{-t}
=+\infty .
$$
In other words, this operator is just not defined. 
This is the point where a coboundary $b$ of \eqref{potential coboundary} shows its significance. 

We recall that 
the derivative of a function $f$ with respect to a  
 Riemannian metric  $ d\sigma =\gamma
\,|dz|$ is given by Formula \eqref{derivative R metric}.
The associated geometric potential is $$\psi = -t\log |f'|_\sigma = -t\log |f'| +t\log \ga -t\log \ga \circ f\,.$$
Since the Euclidean metric plainly does not work, one can try the spherical metric
$d\sigma = |dz|/(1+|z|^2)$ which is another natural choice. Considering again  $f(z)=\l e^z$, we get
$$
\pft \1(w) = \left( \frac{1+|w|^2}{|w|}\right)^t \sum_{f(z)=w}(1+|z|^2)^{-t}
$$
which, this time, is finite provided that $t>1/2$. In fact then, for large $w$, and with $x_0= \log |w/\l |$,
$$
\pft \1(w) \asymp |w|^t \int _\R \frac{dy}{(1+|x_0+iy|^2)^t} =  |w|^t  (1+x_0^2)^{1/2-t }\int _\R \frac{dy}{(1+ y^2)^t} <+\infty,
$$
but 
$$
\lim_{w\to \infty} \pft \1 (w) =+\infty.
$$ 
Thus, $\pft$ is not a bounded operator.

\smallskip

It turns out that for the exponential family and in general for entire functions the logarithmic metric $d\sg = 1/ (1+|z|)$
is best appropriate. This is a natural choice for several reasons. For example, this point of view
is used in Nevanlinna's  value distribution theory. Also, in the dynamics of entire functions from class $\cB$, Eremenko--Lyubich \cite{EL92} have introduced logarithmic coordinates, which now is
a standard tool. Either working in these coordinates or considering derivatives with respect to the logarithmic metric are equivalent things.

\smallskip

\subsection{Growth conditions}
The situation is different for meromorphic functions because of their behavior at poles.
If $f:\C\to\oc$ is meromorphic and if $b$ is a pole of multiplicity $q$, which
 is nothing else than a critical point of multiplicity $q\ge 1$ of $f$, 
 then 
 \beq \label{new 1}
 |f'(z)|\asymp \frac 1{|z-b|^{q+1}} = \frac 1{|z-b|^{q(1+1/q)}}\asymp |f(z)|^{1+1/q} \quad\text{near}\quad b.
 \eeq
Motivated by the exponential family $\lam e^z$, we introduced in \cite{MyUrb08} and \cite{MUmemoirs} some classes of meromorphic functions for which there are relations between $|f'|$ and $|f|$. More precisely:
 
\bdfn[\it Rapid derivative growth and dynamical semi–regularity]
 A meromorphic function $f:\C\to\oc$ is said to have a rapid derivative growth if and only if
 there are $\underline{\al}_2 > \max\{0, -\al_1\}$ and $\kappa>0$ such that
 \begin{equation}\label{eq intro growth}
 |f'(z)| \geq  \ka^{-1}(1+|z|)^{\al _1}(1+|f(z)|^{\underline{\al}_2})
 \end{equation}
for all finite $ z\in \jul \setminus f^{-1}(\infty  )$.
A  \ehyp and expanding meromorphic function of finite order $\rho$ which  satisfies the rapid derivative growth condition is called dynamically semi–regular.
\edfn

 Of course, $f(z)=\l e^z$ satisfies \eqref{eq intro growth}  with $\al_2 \equiv 1$ and $\al_1=0$. 
 The reader can find many other families in Chapter 2 of \cite{MUmemoirs} which are dynamically semi-regular.

 \smallskip
 
 For such functions there is a good choice of the coboundary $b$ or, equivalently, of the 
 Riemannian metric. We recall that we consider metrics of the form \eqref{2.8 full}
 and that we frequently use the simpler form of \eqref{2.8}, namely $d\tau (z) = |z|^{-\tau}|dz|$.
 This is possible as soon as the Fatou set is not empty, which is the case for  \ehyp functions, since then we can assume without loss of generality that $0\in \fat$ and then ignore what happens 
near the origin. 

If $f$ has balanced growth then, setting $\hat \tau=\al_1+\tau$,
\beq\label{2.9}
|z|^{\hat \tau} \preceq |z|^{\hat \tau} |f(z)|^{\underline{\a}_2 -\tau}\preceq |f'(z)|_\tau \preceq
|z|^{\hat \tau} |f(z)|^{\overline{\a}_2 -\tau}\; , \quad z\in \jul \setminus f^{-1}(\infty ),
\eeq
the right hand inequality being true under the weaker condition \eqref{eq intro growth}.
Therefore, for a dynamically semi--regular function $f$, we get 
the estimate
$$
\pft \1 (w) \preceq \frac{1}{|w|^{\underline{\a _2}-\tau}} \sum _{z\in f^{-1}(w)} |z|^{-\hat{\tau} t}
  \quad , \; w\in \jul .
$$
This last sum, which also is called Borel sum, is very well known in Nevanlinna theory.
The order  $\rho$ of $f$ is precisely the critical exponent for this sum. Hence, if $f$ has finite order $\rho$ and  if 
$\hat{\tau} t > \rho$, then it is a convergent series and in fact one has the crucial 
following property (see Proposition 3.6 in \cite{MUmemoirs}).

\bprop\label{prop bdd and decay}
If $f$ is satisfies the rapid derivative growth condition, if $0\in \fat$, and if  $\tau \in (0,\underline{\al}_2 )$ then,  for every $ t > \rho / \hat{\tau}$, there exists $M_t$ such that
\beq\label{13/3 2}
 \pft \1 (w)\leq M_t \quad \text{and}\quad \lim_{w\to\infty} \pft \1 (w) =0
 \quad , \quad w\in \jul \, ,
\eeq
\eprop

Once having property \eqref{13/3 2}, one can develop a full thermodynamic formalism
provided that the function $f$ is  \ehyp and expanding. The first issue is again about
the existence of conformal measures. It is taken care of by Theorem~\ref{t120200730}. Therefore, for  \ehyp and expanding meromorphic functions satisfying the hypotheses of Proposition~\ref{prop bdd and decay}, we have good conformal measures for all $t> \Th$.

We recall that dynamically semi-regular functions have been introduced in Definition \ref{dfn regular fcts}.
The following Perron--Frobenius--Ruelle Theorem is part of Theorem 1.1 in \cite{MyUrb10.2} and
Theorem 5.15 of \cite{MUmemoirs}, which is true for a class of more general potentials.

\bthm \label{theo pfr}
If $f:\amsc\to \cbar$ is a dynamically semi-regular
meromorphic function
then, for every $t>
\frac{\rho }{\hat \tau}$, the following are true.
\begin{itemize}
    \item[(a)] The topological pressure $\P(t)=\lim_{n\to\infty}
    \frac{1}{n} \pft^n(\1)(w)$ exists and is independent of $w\in J(f)$.
    \item[(b)] There exists a unique $\l|f'|_\tau^t$-conformal 
    measure $m_t$ and necessarily $\l=e^{\P(t)}$. 
    
    \item [(c)] There exists
    a unique Gibbs state $\mu_t$ of the parameter $t$, where being Gibbs means that
    $\mu_t$ is a Borel probability $f$--invariant measure absolutely continuous with respect to $m_t$. Moreover, the measures $m_t$ and $\mu_t$
are equivalent and are both ergodic and supported on the conical limit
    set of $f$.
    \item[(d)] The Radon--Nikodym derivative $\psi_t=d\mu_t/dm_t:J(f)\to[0,+\infty)$ is a continuous nowhere vanishing bounded function  satisfying
    $\lim_{z\to\infty}\psi_t(z)=0$.
\end{itemize}
\ethm

Starting from this result, much more can be said but under the stronger growth condition \eqref{eq intro growth beta}. Namely, the Spectral Gap property along with its applications: 

\sp\fr
- The Spectral Gap \cite[Theorem 6.5]{MUmemoirs}

\bthm\lab{t6120101}
If $f$ is a dynamically semi-regular function and if $t>\frac{\rho }{\hat \tau}$, then the following are true.

\,

\begin{itemize}
\item[(a)] The number $1$ is a simple isolated eigenvalue of
the operator $\hat\pf_t:=e^{-\P(t)}\pf_t:\H_\b\to \H_\b$, where $\b\in(0,1]$ is arbitrary and $\H_\b$ is the Banach space of all complex--valued bounded 
H\"older continuous defined on $\jul_f$, equipped with the corresponding H\"older norm. The rest of the spectrum of $\pf_t$ is  contained
in a disk with radius strictly smaller than $1$. In particular, the operator $\hat\pf_t:\H_\b\to \H_\b$ is quasi--compact.
\item[(b)] More precisely: there exists a bounded linear operator $S:\H_\b\to \H_\b$ such that
$$
\hat\pf_t=Q_1+S,
$$
where $Q_1:\H_\b \to \C \den$ is a projector on the eigenspace $\C
\den$, given by the formula 
$$
Q_1(g)=\lt(\int g \, dm_\phi\rt)\den_t,
$$
$Q_1\circ S=S\circ Q_1=0$ and
$$
||S^n||_\b\le C\xi^n
$$
for some constant $C>0$, some constant $\xi\in (0,1)$ and all $n\ge 1$.
\end{itemize}
\ethm

\

\fr
- \cite[Corollary 6.6]{MUmemoirs}

\bcor \label{5.1.2}
With the setting, notation, and hypothesis of Theorem~\ref{t6120101} we have, for every integer $n\geq 1$, that $\npf ^n = Q_1 +S^n$
and that $\npf ^n (g)$ converges to $\left( \int g\, dm_\phi \right) \den $ exponentially fast when $n\to \infty$. More precisely,
$$\lt\| \npf ^n (g) -\left( \int g\, dm_\phi \right) \den \rt\| _\b =\| S^n (g) \| _\b \leq C \xi ^n \|g\| _\b \quad , \;\; g\in H_\b.$$
\ecor

\

\fr
- Exponential Decay of Correlations \cite[Theorem 6.16]{MUmemoirs} 

\bthm\lab{t3120301}
With the setting, notation, and hypothesis of Theorem \ref{t6120101} there exists a large class of functions $\psi_1$ such that for all $\psi_2\in L^1(m_t)$ and all integers $n\ge 1$, we have that
$$
\lt|\int(\psi_1\circ f^n\cdot\psi_2)\,d\mu_t- \int\psi_1\,d\mu_t\int\psi_2\,d\mu_t\rt|\le
O(\xi^n),
$$
where $\xi\in(0,1)$ comes from Theorem~\ref{t6120101}(b), while the big ``O'' constant depends on both $\psi_1$ and $\psi_2$.
\ethm

\

\fr
- Central Limit Theorem \cite[Theorem 6.17]{MUmemoirs}

\bthm
With the setting, notation, and hypothesis of Theorem \ref{t6120101} there exists a large class of functions $\psi$ such that the sequence of random variables
$$
\frac{\sum_{j=0}^{n-1}\psi\circ f^j-n \int\psi\,d\mu_t}{\sqrt{n}}
$$
converges in distribution, with respect to the measure $\mu_t$, to the Gauss (normal) distribution ${\mathcal N}(0,\sigma^2)$ with some $\sg>0$. More precisely,
for every $t\in \R$,
$$
\begin{aligned}
\lim_{n\to\infty}\mu_t\Bigg(\bigg\{z\in \jul_f: \frac{\sum_{j=0}^{n-1}\psi\circ f^j(z)-n \int\psi\,d\mu_t}{\sqrt{n}}&\le t\Bigg\}\Bigg)= \\
&= {1\over \sigma \sqrt{2\pi}}\int_{-\infty}^t\exp\lt(-\frac{u^2}{2\sigma^2}\rt)\, du.
\end{aligned}
$$
\ethm

\

\fr
- Variational Principle \cite[Theorem 6.25]{MUmemoirs}
\bthm \label{thVariational Principle}
With the setting, notation, and hypothesis of Theorem~\ref{t6120101}, we have that
$$
\P(t)=\sup\lt\{{\rm h}_\mu(f)-t\int_{\jul}\log|f'|_1\, d\mu\rt\},
$$
where the supremum is taken over all Borel probability $f$-invariant
ergodic measures $\mu$ with $\int_{\jul}\log|f'|_1\, d\mu>-\infty$. Furthermore, $\int_{\jul}\log|f'|_1\,d\mu_t>-\infty$ and $\mu_t$ is the only one among such measures satisfying the equality
$$
\P(t)={\rm h}_{\mu}(f)-t\int_{\jul}\log|f'|_1\, d\mu.
$$
In the common terminology this means that the $f$--invariant measure $\mu_t$ is the only equilibrium state of the potential $-t\log|f'|_1$.
\ethm

%-------------------------------------------------------------------------------

In  \cite{MUmemoirs} appears also a stronger symmetric growth condition. It is the following 
and it was used in order to get more geometric informations out of the thermodynamical formalism.
The principal application of it was to obtain a Bowen's Formula expressing the hyperbolic dimension
as the zero of the topological pressure function.

\bdfn[\it Balanced growth and dynamical regularity]
A meromorphic function $f:\C\to\oc$ is balanced if and only if there are $\kappa>0$,
a bounded function
$\al_2 : \jul \cap \C \to [\underline{\al}_2 ,
\overline{\al}_2]\subset (0,
\infty)$ and
$\a _1 >-\underline{\al}_2 =-\inf \a_2$
 such that
\begin{equation}\label{eq intro growth beta}
 \ka^{-1}(1+|z|)^{\al _1}(1+|f(z)|^{\al _2(z)})\leq
 |f'(z)| \leq  \ka(1+|z|)^{\al _1}(1+|f(z)|^{\al _2(z)})
\end{equation}
for all finite $ z\in \jul \setminus f^{-1}(\infty  )$.
A balanced  \ehyp and expanding meromorphic function  of finite order $\rho$ is called dynamically regular. \edfn

In this stronger symmetric condition it is important that $\al_2$ is a function since \eqref{new 1}
 shows that  at poles of a meromorphic function this exponent $\al_2$ does depend on the multiplicity $q$. 
 Typical meromorphic functions that satisfy the balanced growth condition
 are all elliptic functions. Again, many other families appear in Chapter 2 of \cite{MUmemoirs}.

%-------------------------------------------------------------------------------

\sp

\subsection{Geometry of tracts} 
For entire functions the thermodynamical formalism is known to hold in a much 
larger setting than the functions that satisfy the growth conditions since we now have a quite optimal approach of \cite{MUpreprint3}. It shows that the geometry of the tracts determines
the behavior of the transfer operator. Let us briefly recall and explain this now.

As it was explained right after the Definition \ref{defi disjoint type},
in order to study the dynamics of a disjoint type entire function $f$ near the Julia set, only its restriction
to the tracts is relevant. Let us here consider the simplest case where $f\in \cB$ 
has only one tract $\Om$. Remember that $f_{|\Om} = e^{\ph^{-1}}$. A simple calculation gives
$$
|f'|_1^{-1}= \frac{|\ph '|}{|\ph |} \circ \ph^{-1} 
$$
in $\Om$. This gives that 
\beq\label{22/4 1}
\pft \1 (w) = \sum_{\xi\in \exp^{-1}(w)} \left|\frac{\ph '(\xi)}{\ph (\xi )}\right|^t
\eeq
 entirely does depend on the conformal representation $\ph$ of the tract
and thus entirely on the tract $\Om$ itself. In fact, the operator $\pft$ does depend
on the geometry of $\Om$ at infinity. 
In order to study the behavior of this operator, one considers the rescaled maps
$$
\ph_T := \frac{1}{|\ph (T)|} \; \ph \circ T : Q_1\longrightarrow \frac{1}{|\ph (T)|}\Om_T
$$
where $Q_T$, especially $Q_1$, has been defined in \eqref{rectangles} and where for $T\geq 1$,
$$
\Om_T:=\ph (Q_T).
$$ 
These maps behave especially well as soon as the tract has some nice geometric properties.

\subsubsection{H\"older tracts} Loosely speaking, a H\"older domain is the image of the unit disk by a H\"older map. But such domains are clearly bounded 
whereas logarithmic tracts are unbounded domains. Following \cite{My-HypDim}, we therefore consider natural exhaustions of the tract by H\"older domains and a scaling invariant
notion of H\"older maps.
A conformal map $h:Q_1\to U$  is called $(H,\al )$--H\"older if and only if
\beq\label{holder 1'}
|h (z_1)-h(z_2)|\leq H |h' (1)| |z_1-z_2|^\al \quad\text{for all}\quad z_1,z_2\in Q_1\,.
\eeq

\bdfn\label{Holder tract}  
 The tract  $\Tract $ is  H\"older, if and only if\eqref{(4.2)} holds and
 the maps $\ph_T$ are uniformly H\"older, i.e.  there exists $(H,\al )$ such that 
 for every $T\geq 1$ the map $\ph_T$ satisfies \eqref{holder 1'}.
\edfn

\fr 
Quasidisks and John domains serve as good examples of H\"older tracts.

\smallskip

\subsubsection{Negative spectrum} The boundary $\partial \Om$ of a tract is an analytic curve. However,
seen from infinity such a boundary may appear quite fractal. In order to quantify this property,
we associate to a tract a version of integral means spectrum (see  \cite{Makarov98} and \cite{PommerenkeBook} for the classical case).
%\begin{figure}[h]
%   \includegraphics[height=2cm]{RescalingsZoom.png}
%     \caption{The part of the boundary in the boxes can resemble more and more a fractal as $T\to\infty$.}
 %    \label{Figure 2}
%\end{figure}
In order to do so, let $h:Q_2\to U$ be a conformal map onto a bounded domain $U$ and define
\beq\label{4}
\b_h(r,t):= \frac{\log \int_I |h'(r+iy)|^tdy}{\log 1/r} \; \text{ ,} \  \ r\in (0, 1) \ \text{ and } \  t\in \R \,.
\eeq
The integral is taken over $I=[-2,-1]\cup[1,2]$ since this corresponds to the part of the boundary of $U$ that is important for our purposes.

Applying this notion to the rescalings $\ph_T$ and then letting $T\to \infty$ leads to desired integral means of the tract $\Om$,
\beq\label{5a}
\b_\infty (t):= \limsup_{T\to +\infty} \b_{\ph_T}(1/T,t),\eeq
and to the associated function
\beq\label{16}
 b_\infty (t):=\b_\infty (t)-t+1 \quad , \quad t\in \R\,.
\eeq
It turns out that the function $b_\infty$ is convex, thus continuous, with 
$
b_\infty (0)=1$ and with $  b_\infty (2)\leq 0$. Consequently, the function $b_\infty$ has at least one zero in $(0,2]$ and we can introduce a number 
$\tstar\in (0,2]$ by the formula
\beq\label{7}
\tstar:= \inf \{ t>0 \; : \; b_\infty (t)=0\}
= \inf \{ t>0 \; : \; b_\infty (t)\le 0\}\,.
\eeq
We only considered here the case of a single tract and the adaption for functions in $\cD$ having finitely many tracts 
is straightforward (see \cite{MUpreprint3}).

\bdfn\label{18} 
A function $f\in \cD$ has negative spectrum if and only if, for every tract,
$$ b_\infty (t)<0 \; \text{ for all} \ \ t> \tstar\,.$$
\edfn

\fr 
A relation of the H\"older tracts property and the negative spectrum  property is provided by the following. 

\bprop[Proposition 5.6 in \cite{MUpreprint3}]\label{28}
A function $f\in\cB$  has negative spectrum if it has only
 finitely many tracts and all these tracts are H\"older.
\eprop

%\smallskip

\subsubsection{Back to the thermodynamic formalism and its applications} From now on we
assume that $f$ is a function of the class $\cD$ and has negative spectrum. Let $\Th$ be again the parameter introduced in \eqref{7}.

Starting from the formula  \eqref{22/4 1}, one can express the transfer
operator  in terms of integral means (see \cite[Proposition 4.3]{MUpreprint3}) in the following way:
\beq\label{pft precise formula}
\pft \1 (w) \asymp (\log |w|)^{1-t} \left\{\int_{-1}^1\left|\ph_{\log|w|} '(1+iy)\right|^t  dy+ 
\sum_{n\geq 1} 2^{n \big(1-t+\b_{\ph_{2^n\log|w|}}( 2^{-n},\, t )\big)}\right\}
\eeq
for every $t\geq 0$ and every $w\in\Om$. The series appearing in this formula may diverge. Nevertheless,
this formula very well describes the behavior of the transfer operator. 
It allows us to develop the thermodynamic formalism if the negative spectrum
assumption holds. The first step is to verify again the conclusion of Proposition \ref{prop bdd and decay}
which, we recall, is crucial for establishing the existence of conformal measures.
Then one can adapt the arguments of \cite{MUmemoirs} to get the following version of the Perron-Frobenius--Ruelle
Theorem (\cite[Theorem 1.2]{MUpreprint3}).

\bthm\label{21/3 1}
Let  $f\in \cD$ be a function having negative spectrum and let $\tstar\in  (0,2]$  be the smallest zero of $b_\infty$. Then, the following hold:
\begin{itemize}
%\addtolength{\itemindent}{-0.3cm}
  \setlength\itemsep{1mm}
\item[-]  For every $t>\tstar$,  the whole thermodynamic formalism, along with its all usual consequences holds: the Perron-Frobenius-Ruelle Theorem, the Spectral Gap property along with its applications: Exponential Mixing, Exponential Decay of Correlations and Central Limit Theorem.
\item[-]  For every $t<\tstar$, the series defining the transfer operator $\pft$ diverges.
\end{itemize}
\ethm

In many cases, by using a standard bounded distortion argument, this result and all its consequences can be extended beyond the class of disjoint type to larger subclasses of hyperbolic functions.
For example, it does hold for all hyperbolic functions in class $\cS$ having finitely many tracts and no necessarily being of disjoint type. This is for example the case for functions of finite order that satisfy \eqref{(4.2)}. The later is a very general kind of 
 quasi--symmetry condition. 
 
\bques
Is the assumption  \eqref{(4.2)} necessary?
\eques

An important feature of the H\"older tract property is is that it is a quasiconformal invariant notion.
This has several important applications. Let us just mention one of them.

\bthm[Theorem 1.3 in \cite{MUpreprint3}] \label{21/3 2} Let $\cM$ be an analytic family of entire functions in class $\cS$. Assume that there
is a function $g\in \cM$ that has finitely many tracts over infinity and that all these tracts are H\"older.
 Then every function $f\in \cM$ has negative spectrum and the thermodynamic formalism holds for every hyperbolic map from $\cM$.
\ethm

\smallskip

Theorem \ref{21/3 1} gives no information at the transition parameter $t=\Th$. For all classical functions
the transfer operator is divergent at $\Th$, and thus the pressure $\P (\Th )=+\infty$. This then implies that the pressure
function has a zero $h>\Th$. Functions with a completely different behavior have been found recently in \cite{MyZd}.

\bthm[\cite{MyZd}]\label{thm myzd1} 
For every $1<\Th<2$ there exists an entire function $f\in\cB$ with the following properties:
\ben 

\item[(a)] The entire function $f$ is of finite order and of disjoint type.

\,
 
\item[(b)] The corresponding transfer operator has transition parameter $\Th$.

\,
 
\item[(c)]  The transfer operator is convergent at $\Th$ and the property \eqref{13/3 2} holds. 

\item[(d)] Consequently,
the Perron--Frobenius--Ruelle Theorem \ref{theo pfr} and its consequences hold at $t=\Th$.

\,
 
\item[(e)] The topological pressure at $t=\Th$ is  strictly negative.

\,
 
\item[(f)] Consequently, the topological pressure of $f$ has no zero. 
\een
\ethm

\fr For the special case of $\Th =2$, the reader can find examples in \cite{Rempe-HypDim2}.

\sp

%{\blue ??}

Here are two more questions related to this section. First of all, we have seen in Proposition \ref{28} that the H\"older tract property
implies negative spectrum. 
For some special functions, Poincar\'e linearizer, both properties coincide (\cite[Theorem 7.8]{MUpreprint3}).
\bques
Are all tracts of any entire function in class $\cB$ with negative spectrum
 H\"older? 
\eques

For H\"older tracts with  corresponding H\"older exponent $\al\in  (1/2, 1]$ it is known that $\Th < 2$. 
\bques
What abut the general case? More precisely, if $\Om$ is a H\"older tract with H\"older exponent $\al\in  (0, 1]$, do we then have that $\Th < 2$? If so,
this would be an analogue of the Jones--Makarov Theorem
\cite{JM95} which states that the Hausdorff dimension of the boundary of an $\al$--H\"older domain is less than two, furthermore,  less than
$2-C\al$ where $C>0$ is a universal constant.
\eques

\section{Hyperbolic dimension and Bowen's Formula}  

The Hausdorff dimension, and in fact all other fractal dimensions, of the Julia set of meromorphic functions have been studied a lot.
The interested reader can consult the survey by Stallard \cite{Stallard-Survey}. Here we focus on the hyperbolic dimension.

\subsection{Estimates for the hyperbolic dimension}
We recall that the \emph{hyperbolic dimension} $\HypDim(f)$ of the function $f$ is the supremum of the Hausdorff dimensions
of all forward invariant compact sets on which the functions is expanding. Right from the definition,
$$\HypDim(f)\leq \HD (\jul ).$$
It has recently been observed by Avila-Lyubich in \cite{AvLyu2} that there are polynomials for which there is strict inequality 
between these two dimension. 
\bthm[\cite{AvLyu2}] There exists a Feigenbaum polynomial $p$ for which 
$$\HypDim(p) < \HD (J(p)) =2.$$
\ethm
\fr Although this result being rather exceptional for rational functions, it appears quite often for transcendental, especially
entire, functions. Stallard \cite{St1990} observed this implicitly and Urba\'nski--Zdunik in \cite{UZ03}.

\bthm[\cite{St1990}, \cite{UZ03}]\label{thm gap hypdim}
There are (even hyperbolic) entire functions $f$ of finite order and of class $\cS$ for which 
$$\HypDim(f) < \HD (\jul ) =2.$$
\ethm
The equality $ \HD (\jul ) =2$ goes back to McMullen's result \cite{McM87}. In either case, of rational functions as well as of
transcendental functions, we do not know any such example with the Hausdorff
dimension of the Julia set equal to $2$. Thus:

\bques
Is there an entire or meromorphic function $f\in\cB$ with with a logarithmic tract over infinity and such that
$$
\HypDim (f)< \HD (\jul )<2 \;\; ?
$$
\eques

While the hyperbolic dimension of a meromorphic functions is often strictly smaller then the dimension of its Julia set. However, it can not be too small as long as 
the function has a logarithmic tract over infinity. In fact 
Bara\'nski, Karpi\'nska and Zdunik \cite{BKZ-2009} obtained the following very general result.

\bthm[\cite{BKZ-2009}]
The hyperbolic dimension of the Julia set of a meromorphic function with a logarithmic tract over infinity is greater than 1. 
\ethm 

For $f_\l (z)=\l e^z$, Karpinska \cite{karpinska-1999} showed that the hyperbolic dimension goes to one as $\l$ goes to zero.
In this sense, the above estimate is sharp. However, if the logarithmic tracts have some regularity then one gets more information,
see \cite{My-HypDim}.

\bthm[\cite{My-HypDim}]\label{thm main intro Holder}
If a meromorphic map $f$ has a logarithmic tract over infinity
and if this tract is H\"older, then
$$\HypDim (f)\geq \Theta \geq 1$$
where $\Th$ is the number defined in \eqref{7}.
\ethm

In this result, one can not expect strict inequality except if $\Th =1$. Indeed,
 for every given $\Th\in (1,2)$ there is an entire function $f$ with H\"older tract such that $\HypDim(f)=\Th$
(see \cite{MyZd}). On the other hand, the paper \cite{My-HypDim} provides a sufficient condition, expressed in terms of
the boundary of the tract, which implies strict inequality. 

\sp

The hyperbolic dimension can also be maximal. This has been shown by Rempe--Guillen \cite{Rempe-HypDim2}.
He first constructs a local version, called now model, and then approximates it by entire functions. His approximation
result is a very precise version of Arakelyan's approximation and is of its own interest.
\bthm[\cite{Rempe-HypDim2}]
There exists a transcendental entire function f of disjoint type and finite order such that $\HypDim(f) = 2$. \\
\ethm

%{\blue On the opposite $.. <2$ sine family etc + break down equality ???}

%\smallskip

\subsection{Bowen's Formula}
The pressure function $t\mapsto\P_\tau(t)$ is convex, hence continuous and, when the map $f$ is expanding, it is also strictly decreasing.
Consequently, there exists a unique zero $h$ of $\P_\tau$ provided 
$$\P_\tau(t)\geq 0 \quad \text{for some $t$.}$$
It goes back to Bowen's paper \cite{Bow79} that this zero is of crucial importance when studying fractal dimensions of limit and Julia sets.
Bowen showed that this number $h$  is the Hausdorff dimension of the limit set for any co-compact quasifuchsian group. His result extends easily to the case of 
of Julia sets of hyperbolic rational functions. 
Since then his formula has been generalized in many various ways
and it became transparent that that for transcendental functions his formula detects the hyperbolic dimension rather than the Hausdorff dimension of the entire Julia sets.

The first result of this kind for transcendental functions is, up to our knowledge, was obtained in \cite{UZ03} and \cite{UZ04} while the 
most general Bowen's Formula for transcendental functions
is due to Bara\'nski, Karpinska and Zdunik \cite{BKZ-Bowen}. Here again, we only formulate a version for  \ehyp functions while their result holds in much bigger generality.

\bthm [\cite{BKZ-Bowen}]
For every  \ehyp meromorphic function $f\in \cB$ we have $\P_{sph} (2)\leq 0$ and 
$$
\HypDim (f) = \HD(\rad ) = \inf \big\{ t>0 \; ;\;\; \P_{sph } (t) \leq 0\big\}\,.
$$
\ethm

We recall that the authors showed the existence of the spherical pressure (Theorem \ref{thm sph pressure})
and that there exists $\Th$ such that the pressure is finite for all $t>\Th$ and infinite for all $t<\Th$. If $\P_{sph } (\Th ) \geq 0$, 
then the pressure has a smallest zero $h\geq \Th$ and this number $h$ turns out to be the hyperbolic dimension. Otherwise,
so if $\P_{sph } (\Th )< 0$, then $\HypDim (f) =\Th$ and in fact such possibility does happen (Theorem \ref{thm myzd1}).

\sp

Other versions of Bowen's formula, with pressure taken with respect to adapted Riemannian metrics, still of the form \eqref{2.8 full},
are contained in \cite{MyUrb08, MUmemoirs, MUpreprint3} and also a version for random dynamics of transcendental functions in \cite{MyUrb2014} and \cite{UZ18}. All these papers contain many other results related to Bowen's formula and formed an important step between \cite{UZ03}, \cite{UZ04} and \cite{BKZ-2009}.

\smallskip

\section{Real Analyticity of Fractal Dimensions} 
Bowen's Formula determines the hyperbolic dimension of a given ``sufficiently hyperbolic'' meromorphic function $f$. But 

\smallskip

\centerline{\emph{what happens to this
dimension when the map $f$ varies in an analytic family?} }

\smallskip

\noindent
For rational functions, this has been explored in detail.
In contrast to the case of entire functions, the radial and Julia sets of a hyperbolic rational function coincide and consequently
also do the corresponding dimensions. Therefore, one is naturally interested in the behavior of the map
$$
f\longmapsto \HD (\jul ).
$$
In 1982, Ruelle \cite{Rue82} positively confirmed a conjecture of Sullivan and showed that the Hausdorff dimension of the Julia set of hyperbolic rational functions depends real--analytically on the map. The hyperbolicity hypothesis is essential here; see \cite[Remark 1.4 ]{Shi98} and also \cite{DSZ97}.

The first result on analytic variation of the hyperbolic dimension of transcendental functions is due to Urba\'nski and Zdunik \cite{UZ03} and concerns the
exponential family $\l e^z$. Since then this property has been obtained for many families
of dynamically regular functions (\cite{MyUrb08, SU14} and \cite{MUmemoirs}; the last of these papers treating also real analyticity 
of appropriate multifractal spectra; for entire functions in class $\cD$ see \cite{MUpreprint3}). For the same kind of families,
such analyticity is also true in the realm of random dynamics; see \cite{Mayer2016aa}.

Instead of presenting a complete overview of all relevant, sometimes quite technical, results we now describe the general framework
followed by two representative methods and results.

\smallskip

Similarly as the hyperbolicity hypothesis
for rational functions, there are  a number of conditions, in a sense necessary,
needed to expect real analytic variation of the hyperbolic dimension in the transcendental case. They can be summarized as follows. 
\ben
\item[-] $\cF$ is an analytic family of meromorphic functions. The reader simply can assume that 
$\cF = \{f_\l=\l f\; : \;\; \l\in \La\}$ where $f$ is a given meromorphic function and $\La$ an open subset of $\C^*$.
%In order to keep the presentation of this text less technical  we allow us to simply consider
 %families of the type $\cF = \{\l f\; , \;\; \l\in \La\}$ where $f$ is a given meromorphic function and $\La$ an open subset of $\C^*$.
 Clearly there are more general settings. For example, in the case of entire functions in class $\cS$ there is a natural notion of 
analytic family due to Eremenko--Lyubich \cite{EL92}; they are in particular always  finite dimensional. 
\item[-] The functions of $\cF$ are   \ehyp and expanding.
\item[-] The family $\cF$ is structurally stable in the sense of holomorphic motions. 
\een

In most results the holomorphic motion is also assumed to have some uniform behavior which is for example implied by a condition called  bounded deformation \cite{MUmemoirs}. 

The last commonly used hypothesis is that
the full thermodynamic formalism applies. Here appears a crucial fact which is specific to the transcendental case. 
%Contrary to the case of a rational functions, 
The transfer operator of a transcendental function is usually not defined for small parameters $t>0$.
Let us follow the notation used in Theorem \ref{21/3 1}
and call again $\Th$ the transition parameter. In fact, one must rather write $\Th_f$ since this number can depend on a particular function $f$
from  a given family $\cF$.
\ben
\item[-] The thermodynamic formalism holds for the functions in  $\cF$ with constant transition  parameter
$\Th = \Th _{f_\l }$, $ \l \in \La.$ 
\een

In particular, Bowen's Formula applies to the functions we consider here and thus two cases appear: for $f\in \cF$, either
\beq\label{cases theta}
\HypDim (f) > \Th \quad \text{ or }\quad \HypDim (f) = \Th .
\eeq

\sp
The first analyticity result we present here is due to Skorulski--Urba\'nski obtained in \cite{SU14}.

\bthm[\cite{SU14}]\label{theo SkUrb}
Suppose that $\Lambda\subset \C$ is an open set, $\cF=\{f_\l \}_{\l\in\Lambda}$ is an analytic family  of meromorphic functions
and that, for some $\l_0\in \Lambda$,
$f_{\l_0}: \C\to \cbar$ is a dynamically regular meromorphic function with  $\HypDim \lt(f_{\lam_0}\rt) > \Th_{f_{\lam_0}}$ and which belongs to class S. Then the function 
$$
\l\longmapsto \HD(J_r(f_\l ))
$$ 
is real--analytic in some open neighborhood of $\l_0$.
\ethm

%{\blue Mariusz: please adapt if not as you like:}

Notice that here the main hypotheses are only imposed on the function $f_{\l_0}$ and not on all functions 
in a neighborhood of it. The authors obtained this result by associating to the globally defined functions
locally defined iterated functions systems (IFS). This is possible by employing so called
nice sets whose existence in the transcendental case is due to Doobs \cite{Do11} and which have been  initially
 brought to complex
dynamics by Rivera--Letelier in \cite{Ri} and Przytycki and Rivera--Letelier in \cite{PrRL07}. An open connected set $U\subset\C$ is called nice if and only if every connected component of $f^{-n}(U)$ is
either contained in $U$ or disjoint from $U$. If $U$ is disjoint from the post--singular set,
then one can consider all possible holomorphic inverse branches of iterates of $f$ and the properties of the nice set imply that
the inverse branches that land in $U$ for the first time define a good countable alphabet conformal IFS in the sense of \cite{MU1} and\cite{MauUrb03}.
It turns out that the limit set of this IFS has the same dimension as the
hyperbolic dimension of $f$ \cite[Theorem 3.4]{SU14}. Thus it suffices to consider IFSs. The later have been extensively 
studied \cite{MauUrb03} providing many useful tools, and, especially, developing the full thermodynamic formalism, and introducing the concepts of regular, strongly regular, co-finitely regular and irregular conformal IFSs. One of the greatest challenges
to apply Theorem \ref{theo SkUrb} is to show that $\HypDim \lt(f_{\lam_0}\rt) > \Th_{f_{\lam_0}}$. In terms of the associated conformal IFSs this means that the IFS coming from $f_{\lam_0}$ is strongly regular.

\smallskip
The common underlying strategy for establishing real analytic variation of Hausdorff dimension of limit sets of conformal IFSs, see \cite{UZ03, MyUrb08, MUmemoirs, MUpreprint3} for ex., is to complexify the setting and to
apply Kato--Rellich Perturbation Theorem. The later is possible thank's to the spectral gap
property which means that $\exp ({\P(t)})$ is a leading isolated simple eigenvalue of the transfer operator and the rest of the spectrum of this operator is contained in a disk centered at $0$ whose radius is strictly smaller than $\exp ({\P(t)})$. An alternative powerful strategy is used in
\cite{Mayer2016aa}. It is based on Birkhoff's approach \cite{Bir57} to the Perron--Frobenius Theorem via positive cones. 
This method has been successfully applied in various contexts. The paper \cite{Mayer2016aa} which deals with random dynamics,
is based on ideas from Rugh's paper \cite{Rug08} who used complexified cones. 
This powerful method works well as soon as appropriate invariant cones are found
and strict contraction of the transfer operator in the appropriate Hilbert metric has been shown. The following is a particular result in \cite{Mayer2016aa}.

\bthm[\cite{Mayer2016aa}]\label{thm intro 1 random}
Let $f_\eta (z)=\eta e^z$  and let $a\in (\frac{1}{3e}, \frac{2}{3e})$
and $0<r< r_{max}$, $r_{max}>0$. Suppose that $\eta_1,\eta_2,..$ are i.i.d. random variables uniformly distributed in $\D (a,r)$. Let $J_{\eta_1,\eta_2,...}$ denote the Julia set of the sequence of compositions 
$$
f_{\eta_n}\circ f_{\eta_{n-1}}\circ\ldots\circ f_{\eta_2}\circ f_{\eta_1}:\C\longrightarrow\C,  \  \  \  n\geq 1,
$$
and let 
$$
J_r (\eta_1,\eta_2,...)= \big\{z\in J_{\eta_1,\eta_2,...}  : \; \liminf_{n\to\infty} |f_{\eta_n}\circ ...\circ f_{\eta_1}(z)|<+\infty\big\}
$$
be the radial Julia set of $\{f_{\eta_n}\circ ...\circ f_{\eta_1}\}_{n\geq1}$. 
Then, the Hausdorff dimension of $J_r(\eta_1,\eta_2,...)$ is almost surely constant and depends real-analytically on the parameters $(a,r)$ provided that $r_{max}$ is sufficiently small.
\ethm

\sp

In contrast to the case of hyperbolic rational functions, analytic variation of the hyperbolic dimension can fail
in the class of hyperbolic entire functions of bounded type. 
This has been recently proved in \cite{MyZd}.

\bthm[\cite{MyZd}]\label{theo main MyZd}
There exists a holomorphic family $\cF =\{f_\l= \l \, f \; , \;  \l\in \C ^*\}$
of finite order entire functions in class $\cB$ 
such that 
the functions $f_\l$, $\l\in (0, 1]$, are all in the same hyperbolic component of the parameter space but the function
$$\l \mapsto {\rm HypDim}(f_\l )$$
is not  analytic in $(0, 1]$.
\ethm

In order to obtain this result, the authors exploited the dichotomy  of \eqref{cases theta}.
In fact, all positive  analyticity results use, sometimes implicitly like in Theorem \ref{thm intro 1 random},
the assumption $\HypDim (f) > \Th$. Using the formula \eqref{pft precise formula} for the transfer operator,
Mayer and Zdunik where able to construct in \cite{MyZd} entire functions for which $\HypDim (f) =\Th$, so obtaining the very special case of equality in \eqref{cases theta}.  Moreover, among these functions there are some that have     strictly negative pressure at $\Th$
which is the key point not only for Theorem \ref{theo main MyZd} but also for the absence of recurrent conformal
measures in Theorem \ref{theo main' MyZd}.

%\bques
%It would be interesting to know if Theorem \ref{theo main MyZd} also holds inside the entire functions of class $\cS$.
%\eques

\smallskip

%\vfill
%\pagebreak

\section{Beyond hyperbolicity}

For many kinds of non--hyperbolic holomorphic/conformal dynamical systems various forms of thermodynamical formalism have been also successfully developed and usually much earlier than for transcendental dynamics. This is the case for rational functions and generalized polynomial--like mappings
having certain type of critical points in the Julia set so that the functions are no longer hyperbolic but
sufficient expansion is maintained. Most notably this is so for parabolic rational functions, subexpanding rational functions, and most generally, for non--recurrent rational functions and topological Collet--Eckmann rational functions; see ex. cite{GPS90}, \cite{P90}, \cite{DU2}, \cite{DU3}, \cite{DU4}, \cite{DU LMS}, \cite{DMNU}, \cite{ADU93}, \cite{DU Subexp}, \cite{DU ACIM}, \cite{DU ES}, \cite{U94}, \cite{U97}, \cite{P98}, \cite{PrRL07}, \cite{PrRL11}, \cite{AvLyu2}, \cite{AL08}, \cite{P18}, \cite{SU03}, \cite{SU04}, and the references therein. Note that some of these papers such as \cite{P90}, \cite{DU ES} and 
\cite{P18} for ex. deal with all rational functions, in particular with no restrictions on critical points at all.

\sp But there is a substantial difference with the hyperbolic case. Except perhaps \cite{DU ES} and \cite{P90}, the Perron--Frobenius (transfer) operator for the original system is then virtually of no use - no change of Riemaniann metric seems to work. The most relevant questions are then about the structure of conformal measures, most notably, their existence, uniqueness, and atomlessness, and about Borel probability invariant measures absolutely continuous with respect such conformal measures, their existence, uniqueness and stochastic properties. Also, application of such results to study the fractal structure of Julia sets. 

\sp Similarly as for non--expanding rational functions, also for non--hyperbolic non--expanding transcendental, entire and meromorphic, functions some forms of thermodynamic formalism have been developed. For the papers coping with critical points in the Julia sets, which is closest to rational functions, see for ex. \cite{KU03}, \cite{KS08b}. One class of trancendental meromorphic functions deserves here special attention. These are elliptic (doubly periodic) meromorphic functions. The first fully developed account of thermodynamic formalism for all elliptic functions and H\"older continuous potentials (satisfying some additional natural hypotheses) was presented in 
\cite{MyUrb05.2}. Up to our best knowledge all other contributions to thermodynamic formalism for elliptic functions deal with geometric potentials of the form $-t\log|f'|$. We would like to mention in this context the paper \cite{KU04}, and, especially, the book \cite{KU20}, which provides an extensive and fairly complete account of thermodynamic formalism for many special, but quite large, classes of elliptic functions with some sufficiently strong expanding features.

\sp The main difficulty and main point of interest in the classes of meromorphic functions discussed in the last paragraph were caused by critical points lying in the Julia sets. Going beyond critical points, there are visible two directions of research. Both of them deal with transcendental entire functions where there are logarithmic singularities, in the form of asymptotic values, in the Julia sets. 

One of them was initiated in \cite{UZ07} dealing with exponential functions $\lambda e^z$, where $0$, the asymptotic value, was assumed to escape to infinity sufficiently fast. The existence and uniqueness of conformal measures and the existence and uniqueness of Borel probability invariant measures absolutely continuous with respect to those conformal measures were proved therein. Its follow up was the paper \cite{UZ18} dealing with analogous classes of functions but iterated randomly. The full (random) thermodynamic formalism with respect to random conformal and invariant measures was laid down and developed therein. 

The second direction of research initiated and developed in \cite{MyUrb10.2} 
aimed to analyze the contribution of non--recurrent logarithmic singularities.
Indeed, the paper \cite{MyUrb10.2} by
 Mayer--Urba\'nski considers the class  of meromorphic functions with polynomial Schwarzian derivatives.
 For example the  tangent family belongs to this class and in general
such functions have no  critical points and they have only finitely many logarithmic singularities. A surprising outcome of this paper was that the behavior of invariant measures absolutely continuous with respect to conformal measures did depend on the order of the function.

\bthm \label{thm main}
Let $f:\oc\to\oc$ be a meromorphic function $f$ of polynomial Schwarzian derivative and assume that it is semi-hyperbolic
in the following sense:
\begin{enumerate}
  \item[-] All the asymptotic values are finite.
  \item[-] The asymptotic values that belong to the Fatou set belong to attracting components.
   \item[-] The asymptotic values that belong to the Julia set have bounded and non-recurrent forward orbits.
   \end{enumerate}
Let $h:=\HD(J(f))$.

   Then,  a Patterson--Sullivan typ construction provides an atomless $h$--conformal measure
   and this  measure is weakly metrically exact, hence
ergodic and conservative. Moreover, there exists a $\sg$--finite invariant measure $\mu$
absolutely continuous with respect to $m$ and this measure
$$
\mu \;\text{ is finite}\quad \text{ if and only if} \qquad h>3\frac{\rho}{\rho +1}
$$
 where $\rho = \rho (f)$ is the order of the function $f$. If $\mu$ is finite, then the
dynamical systems ($f,\mu$) it generates is metrically exact and, in consequence, its Rokhlin's
natural extension is K-mixing.
\ethm

Notice that $3\frac{\rho}{\rho +1}\geq 2$ if and only if the order
$\rho \geq 2$. Consequently the measure $\mu$ is most often infinite. However,
in the case of the
tangent family, which is just one specific example among others, this invariant
measure can be finite.

\vfill
\pagebreak

\bibliographystyle{amsalpha}

%---------------------------------------------------------------------------------------------------------------------------------------------------%

%\backmatter
%\bibliographystyle{spmpsci}
%\bibliographystyle{plain}
%\bibliography{/Users/VM/SynchroCloud/Math/biblio_VM}
%\printindex

%\begin{thebibliography}{10}

%\end{thebibliography}

%---------------------------------------------------------------------------------------------------------------------------------------------------%

\end{document}